\input amstex\documentstyle{amsppt} 
\pagewidth{12.5cm}\pageheight{19cm}\magnification\magstep1
\topmatter
\title $\bold Z/m$-graded Lie algebras and perverse sheaves, II\endtitle
\author George Lusztig and Zhiwei Yun\endauthor
\address{Department of Mathematics, M.I.T., Cambridge, MA 02139,
Department of Mathematics, Stanford University, Stanford, CA 94305}\endaddress
\thanks{G.L. supported by NSF grant DMS-1566618; Z.Y. supported by NSF 
grant DMS-1302071 and the Packard Foundation.}\endthanks
\endtopmatter   
\document

       \define\tfm{\ti{\fm}}
       \define\tfl{\ti{\fl}}

     \define\bco{\bar{\co}}

\define\mpb{\medpagebreak}

 \define\hfp{\hat{\fp}} \define\hfl{\hat{\fl}} \define\hfu{\hat{\fu}}

    \define\hP{\hat P}

\define\frl{\forall}

\define\si{\sim}
\define\wt{\widetilde}
\define\sqc{\sqcup}

\define\ovsc{\overset\cir\to}
\define\qua{\quad}

\define\hL{\hat L}

\define\tcl{\ti\cl}

\define\baf{\bar f}

\define\lb{\linebreak}

\define\op{\oplus}
   
\define\part{\partial}
\define\emp{\emptyset}
\define\imp{\implies}
\define\ra{\rangle}
\define\n{\notin}
\define\iy{\infty}
\define\m{\mapsto}
\define\do{\dots}
\define\la{\langle}
\define\bsl{\backslash}

\define\Lra{\Leftrightarrow}

\define\sub{\subset}    

\define\T{\times}
\define\ti{\tilde}
\define\nl{\newline}
\redefine\i{^{-1}}
\define\fra{\frac}
\define\un{\underline}
\define\ov{\overline}
\define\ot{\otimes}
\define\bbq{\bar{\QQ}_l}

\define\ad{\text{\rm ad}}
\define\Ad{\text{\rm Ad}}
\define\Hom{\text{\rm Hom}}

\define\Ind{\text{\rm Ind}}    \define\tInd{\wt{\Ind}}

\define\ind{\text{\rm ind}}

\define\he{\heartsuit}

\define\a{\alpha}
\redefine\b{\beta}

\define\g{\gamma}
\redefine\d{\delta}
\define\e{\epsilon}
\define\et{\eta}
\define\io{\iota}

\define\ph{\phi}
\define\ps{\psi}
\define\r{\rho}
\define\s{\sigma}
\redefine\t{\tau}
\define\th{\theta}
\define\k{\kappa}
\redefine\l{\lambda}
\define\z{\zeta}
\define\x{\xi}

\define\vp{\varpi}
\define\vt{\vartheta}

\redefine\G{\Gamma}

\define\Si{\Sigma}

\redefine\L{\Lambda}

\define\boc{\bold c}

\define\kk{\bold k}
\define\mm{\bold m}

\define\BB{\bold B}

\define\DD{\bold D}
\define\EE{\bold E}
\define\FF{\bold F}

\define\NN{\bold N}

\define\QQ{\bold Q}

\define\VV{\bold V}

\define\ZZ{\bold Z}

\define\ca{\Cal A}
\define\cb{\Cal B}
\define\cc{\Cal C}
\define\cd{\Cal D}
\define\ce{\Cal E}

\define\ch{\Cal H}
\define\ci{\Cal I}

\define\ck{\Cal K}
\define\cl{\Cal L}
\define\cm{\Cal M}

\define\co{\Cal O}

\define\cq{\Cal Q}
\define\car{\Cal R}
\define\cs{\Cal S}
\define\ct{\Cal T}

\define\cw{\Cal W}
\define\cz{\Cal Z}
\define\cx{\Cal X}

\define\fg{\frak g}
\define\fh{\frak h}

\define\fl{\frak l}
\define\fm{\frak m}

\define\fp{\frak p}
\define\fq{\frak q}

\define\fs{\frak s}

\define\fu{\frak u}

\define\fz{\frak z}

\define\fB{\frak B}

\define\fH{\frak H}

\define\fR{\frak R}

\define\fT{\frak T}

\define\tA{\ti A}

\define\tC{\ti C}

\define\tF{\ti F}

\define\tI{\ti I}

\define\tL{\ti L}

\define\tP{\ti P}

\define\tT{\ti T}

\define\sha{\sharp}

\define\cir{\circ}

\define\GRA{L4}
\define\CUS{L6}
\define\CSV{L7}
\define\CB{L8}
\define\KLL{KL2}
\define\LY{LY}
\define\da{\delta}
\define\doet{\dot\eta}
\head Contents\endhead
10. The vector space $\VV$ and the sesquilinear form $(:)$.

11. The $\ca$-lattice $\VV_\ca$.

12. Purity properties.

13. An inner product.

14. Odd vanishing.

\head Introduction\endhead
As in \cite{\LY} we fix $G$, a semisimple simply connected algebraic group 
over $\kk$ (an algebraically closed field of characteristic $p\ge0$) and a 
$\ZZ/m$-grading $\fg=\op_{i\in\ZZ/m}\fg_i$ for the Lie algebra $\fg$ of $G$; 
here $m$ is an integer $>0$ and $p$ is $0$ or a large prime number. We fix
$\et\in\ZZ-\{0\}$ and let $\da$ be the image of $\et$ in $\ZZ/m$. Let $G_{\un0}$ 
be the
closed connected subgroup of $G$ with Lie algebra $\fg_{\un0}$; this group
acts naturally on $\fg_\da$ and on $\fg_\da^{nil}=\fg_\da\cap\fg^{nil}$ where
$\fg^{nil}$ is the variety of nilpotent elements in $\fg$. 
Recall that $\ci(\fg_{\da})$ denotes the set of pairs $(\co,\cl)$ where $\co$ is a $G_{\un0}$-orbit on $\fg_{\da}^{nil}$ and $\cl$ is an irreducible $G_{\un0}$-equivariant local system on $\co$ up to isomorphism. Let $\fB$ be the set of isomorphism 
classes of simple $G_{\un0}$-equivariant perverse sheaves on $\fg_\da^{nil}$. 
Then there is a natural bijection $\ci(\fg_{\da})\to \fB$ given by $(\co,\cl)\mapsto \cl^{\sha}[\dim\co]$ (for the notation $\cl^{\sha}$ see \cite{\LY, 0.11}).
In \cite{\LY, Theorem 0.6} we have shown that $\fB$ can be naturally decomposed as a 
disjoint union of blocks ${}^\x\fB$ 
indexed by the $G_{\un0}$-conjugacy classes of admissible systems $\x=(M,M_{0},\fm,\fm_{*},\tC)$ 
such that if $B,B'$ are in different blocks
then the $G_{\un0}$-equivariant $\text{\rm Ext}$-groups of $B$ with $B'$ are 
zero. This generalizes a result of \cite{\GRA} in the $\ZZ$-graded case; its
analogue in the ungraded case is the known partition of the set 
$G$-equivariant simple perverse sheaves on $\fg^{nil}$ into blocks given by 
the generalized Springer correspondence.

In this paper we give an (essentially) combinatorial way to parametrize the 
objects in a fixed block ${}^\x\fB$ of $\fB$.
It is known that in the ungraded case, the objects in a fixed block are 
indexed by the irreducible representations of a certain (relative) Weyl group
attached to the block. In the $\ZZ/m$-graded case we will associate to
the block ${}^\x\fB$ a $\QQ$-vector space $\EE$ with a certain finite
collection of hyperplanes. The complement of the union of these hyperplanes
is naturally a union of finitely many (not necessarily conical) chambers 
which can be taken to form a basis of a $\QQ(v)$-vector space $\VV'$ (here $v$
is an indeterminate). 
The chambers represent various spiral induction
associated to the block. The vector space $\VV'$ carries a natural, explicit,
sesquilinear form $(:):\VV'\T\VV'@>>>\QQ(v)$ defined in terms of 
dimensions of $\text{\rm Ext}$-groups 
between the spiral inductions (see \cite{\LY, 6.4}) which correspond to
the chambers. We show that the left radical of this form is the same as its
right radical. Taking the quotient of $\VV'$ by the left or right radical
we obtain a $\QQ(v)$-vector space $\VV$ with an induced sesquilinear form
$(:)$. It turns out that $\VV$ is naturally isomorphic to the Grothendieck 
group based on ${}^\x\fB$ (tensored with $\QQ(v)$) and, in particular, the 
number of elements in ${}^\x\fB$ is equal to $\dim_{\QQ(v)}\VV$.
The vector space $\VV$ has a natural bar involution $\bar{}:\VV@>>>\VV$ and a
natural $\ca$-lattice $\VV_\ca$ in $\VV$, both defined combinatorially. 
(Here $\ca=\ZZ[v,v\i]$.)
We then define a subset $\BB'$ of $\VV$ as the set of all $b\in\VV_\ca$ such
that $\bar b=b$ and $(b:b)\in1+v\ZZ[[v]]$. We show that $\BB'$ is a
signed basis of $\VV$. (Although the definition of $\BB'$ is combinatorial, the
proof of the fact that it is a well defined signed basis is based on geometry;
it is not combinatorial. It would be desirable to find a proof without using
geometry.) Let $\BB'/\pm$ be the set of orbits of the $\ZZ/2$-action $b\m-b$
on $\BB'$. We show that $\BB'/\pm$ is in natural bijection with the given
block ${}^\x\fB$. Thus $\BB'/\pm$ could be regarded as a combinatorial
index set for ${}^\x\fB$. (A similar result in the $\ZZ$-graded case appears
in \cite{\GRA}.)

\mpb

We now discuss the contents of the various sections.
In Section 10 we define the $\QQ$-vector vector $\EE$ with its hyperplane
arrangement associated to a block. We also define the $\QQ(v)$-vector space
$\VV'$ with its sesquilinear form $(:)$, its quotient space $\VV$ and the bar 
operator. In Section 11 we define the $\ca$-lattice $\VV_\ca$ in $\VV$ and the
signed basis $\BB'$. In Section 12 we prove some purity properties of the 
cohomology sheaves of the simple $G_{\un0}$-equivariant perverse sheaves on 
$\fg_\da^{nil}$, which generalize those in the $\ZZ$-graded case given in 
\cite{\GRA}. In Section 13 we generalize an argument in \cite{\GRA} to
express the matrix whose entries are the values of the $(:)$-pairing at two
elements of $\fB$ as a product of three matrices. This is used in Section
14 to prove the vanishing of the odd cohomology sheaves of the intersection 
cohomology of the closure of any $G_{\un0}$-orbit in $\fg_\da^{nil}$ with 
coefficients in an irreducible $G_{\un0}$-equivariant 
local system on that orbit. This generalizes a result in \cite{\GRA} in the
$\ZZ$-graded case whose proof was quite different (it was based on geometric 
arguments which are not available in the present case). 

\mpb

We will adhere to the notation and assumptions of \cite{\LY}. 

\head 10. The vector space $\VV$ and the sesquilinear form $(:)$\endhead
In this section, we introduce a combinatorial way of calculating the number of 
irreducible perverse sheaves in a block ${}^\x\cq(\fg_\da)$. This is achieved 
by introducing a finite-dimensional vector space $\VV'$ over $\QQ(v)$ 
together with a sesquilinear form on it coming from the pairing between spiral 
inductions in the fixed block $\x$. Then the number of irreducible perverse 
sheaves in ${}^\x\cq(\fg_\da)$ turns out to be the rank of this sesquilinear
form on $\VV'$ (Proposition 10.19).

\subhead 10.1\endsubhead
We fix $\x\in\un\fT_\et$ and a representative 
$\dot\x=(M,M_0,\fm,\fm_*,\tC)\in\fT_\et$ for $\x$. We also fix 
$\ph=(e,h,f)\in J^M$ such that $e\in\ovsc\fm_\et$, $h\in\fm_0$, $f\in\fm_{-\et}$.
We set $\io=\io_\ph\in Y_M$. Let 
$$Z=\cz^0_M.$$ 
Since $Z$ is a torus, $Y_Z$ is naturally a free abelian group (with operation
written as addition) and $\EE:=Y_{Z,\QQ}$ may be naturally identified with the
$\QQ$-vector space $\QQ\ot Y_Z$. Let $X_Z=\Hom(Y_Z,\ZZ)$ and let $\la:\ra$ be 
the obvious perfect bilinear pairing $Y_Z\T X_Z@>>>\ZZ$. This extends to a 
bilinear pairing $\EE\T(\QQ\ot X_Z)@>>>\QQ$ denoted again by $\la:\ra$. For 
any $\a\in X_Z$ let
$$\fg^\a=\{x\in\fg;\Ad(z)x=\a(z)x\qua\frl z\in Z\}.$$
For any $(\a,i)\in X_Z\T\ZZ/m$ let
$$\fg^\a_i=\fg^\a\cap\fg_i.$$
For any $(\a,n,i)\in X_Z\T\ZZ\T\ZZ/m$ let
$$\fg^{\a,n}_i=\fg^\a_i\cap({}^\io_n\fg_i).$$
For any $i\in\ZZ/m$, let $\car_i=\{(\a,n)\in X_Z\T\ZZ;\fg^{\a,n}_i\ne0\}$,
$\car^*_i=\{(\a,n)\in\car_i;\a\ne0\}$. Note that 

(a) {\it $\dim\fg^{\a,n}_i=\dim\fg^{-\a,-n}_{-i}$ for any $\a,n,i$; hence
$(\a,n)\m(-\a,-n)$ is a bijection $\car_i@>\si>>\car_{-i}$ and a 
bijection $\car^*_i@>\si>>\car^*_{-i}$.}
\nl
We have 
$$\fg=\op_{(\a,n,i)\in X_Z\T\ZZ\T\ZZ/m}(\fg^{\a,n}_i)=
\op_{i\in\ZZ/m,(\a,n)\in\car_i}(\fg^{\a,n}_i).$$
For any $N\in\ZZ$ and any $(\a,n)\in\car^*_{\un N}$ we set
$$\fH_{\a,n,N}=\{\vp\in\EE;\la\vp:\a\ra=2N/\et-n\}.$$
This is an affine hyperplane in $\EE$. We set
$$\EE'=\EE-\cup_{N\in\ZZ,(\a,n)\in\car^*_{\un N}}\fH_{\a,n,N}.$$

\subhead 10.2\endsubhead
Recall the notation $Y_{H,\QQ}$ for a connected algebraic group $H$ with Lie algebra $\fh$ 
from \cite{\LY, 0.11}. For 
$\mu,\mu'\in Y_{H,\QQ}$, we say $\mu$ commutes with $\mu'$ if for some $r,r'$ in $\ZZ_{>0}$ such that 
$r\mu,r'\mu'\in Y_H$, the images of the homomorphisms $r\mu, r'\mu':\kk^*\to H$ commute with each other. 
This property is independent of the choice of $r,r'$. If $\mu$ commutes with $\mu'$, and 
$r,r'\in\ZZ_{>0}$ are as above, then letting $\l=r\mu, \l'=\r'\mu'$, we have the homomorphism 
$\nu:\kk^*\to H$ given by $\nu(t)=\l(t^{r'})\l'(t^r)$. We then define 
$\mu+\mu':=\nu/(rr')\in Y_{H,\QQ}$. Then $\mu+\mu'$ is independent of the choice of $r,r'$. 
Moreover, 
for any $\k\in\QQ$, we have 
$${}^{\mu+\mu'}_\k\fh=\op_{s,t\in\QQ;s+t=\k}({}^\mu_s\fh\cap{}^{\mu'}_t\fh).
\tag a$$ 
Let $\vp\in\EE$. 
Since $\vp\in\EE=Y_{Z,\QQ}$ commutes with $\io\in Y_{M_0}$, by the above discussion, we may define
$$\un\vp=\fra{|\et|}{2}(\vp+\io)\in Y_{M_{0},\QQ}\sub Y_{G_{\un0},\QQ}.$$
From the definitions, $\un\vp$ may be computed as follows. We choose $f\in\ZZ_{>0}$ such that 
$\l':=f\vp\in Y_Z$ and $f/|\et|\in\ZZ$; we define $\l\in Y_{G_{\un0}}$ by
$\l(t)=\io(t^f)\l'(t)=\l'(t)\io(t^f)$ for all $t\in\kk^*$; we then have
$\un\vp=\fra{|\et|}{2f}\l\in Y_{G_{\un0},\QQ}$. 

We shall set $\e=\doet\in\{1,-1\}$ (the sign of $\et$, see 0.12).

Now ${}^\e\fp^{\un\vp}_*$, ${}^\e\tfl^{\un\vp}_*$ are well-defined, see 2.5, 2.6. 

We set
$${}^\e\tfl^{\un\vp}=\op_{N\in\ZZ}{}^\e\tfl^{\un\vp}_N.$$
We show:

(b) {\it We have ${}^\e\tfl^{\un\vp}_*=\fm_*$ if and only if $\vp\in\EE'$.}
\nl
If $x\in\fg^{\a,n}_i$ and $t\in\kk^*$, then
$$\Ad(\l(t))x=\Ad(\io(t^f))\Ad(\l'(t))x=t^{nf+\la\l':\a\ra}x$$
and $\Ad(\io(t))x=t^nx$. Thus, $\fg^{\a,n}_i\sub{}^\l_{nf+\la\l':\a\ra}\fg_i$ 
and $\fg^{\a,n}_i\sub{}^\io_n\fg_i$. It follows that for any $k\in\ZZ$ we have
$${}^\l_k\fg_i=\op_{(\a,n)\in\car_i;nf+\la\l':\a\ra=k}(\fg^{\a,n}_i),\tag c$$ 
$${}^\io_k\fg_i=\op_{(\a,n)\in\car_i;n=k}(\fg^{\a,n}_i).\tag d$$ 
For any $N\in\ZZ$ we have
${}^\e\tfl^{\un\vp}_N={}^\l_{2fN/\et}\fg_{\un N}$.
We see that
$${}^\e\tfl^{\un\vp}_N=\op_{(\a,n)\in\car_{\un N};nf+\la\l':\a\ra=2fN/\et}
(\fg^{\a,n}_{\un N}).\tag e$$ 
Recall from 1.2(e) that if $N\in\ZZ,\fm_N\ne0$ then $N\in\et\ZZ$. Let $N\in\et\ZZ$. 
From the arguments in 3.3 (or by 10.3(a)) we see that, for some $\vp'\in\EE$, 
we have
${}^\e\tfl^{\un\vp'}_*=\fm_*$. (Here $\un\vp'$ is attached to $\vp'$ in the
same way that $\un\vp$ was attached to $\vp$.) Hence, using (e) with $\vp$ 
replaced by $\vp'$, we see that there exists a subset $R_N$ of $\car_{\un N}$ 
such that $\fm_N=\op_{(\a,n)\in R_N}(\fg^{\a,n}_{\un N})$. If $(\a,n)\in R_N$, 
then $0\ne\fg^{\a,n}_{\un N}\sub\fm_N$. Hence for 
$x\in\fg^{\a,n}_{\un N}-\{0\}$ and $z\in Z$, $t\in\kk^*$,  we have 
$\Ad(z)x=\a(z)x$, $\Ad(\io(t))x=t^nx$; but $\Ad(z)x=x$ since $Z$ is in the 
centre of $M$ and $\Ad(\io(t))x=t^{2N/\et}x$ since $x\in\fm_N$. Thus $\a(z)=1$ 
(so that $\a=0$) and $n=2N/\et$. We see that $R_N\sub\{(0,2N/\et)\}$.
 Conversely, 
we show that $\fg^{0,2N/\et}_{\un N}\sub\fm_N$. If $x\in\fg^{0,2N/\et}_{\un N}$, 
then $\Ad(z)x=x$ for all $z\in Z$ and $\Ad(\io(t))x=t^{2N/\et}x$ for all 
$t\in\kk^*$. Since $\dot\x$ in 9.1 is admissible, there exists $t_0\in\kk^*$, 
$z\in Z$, such that $\fm$ is the fixed point set of 
$\Ad(\io(t_0)z)\th:\fg@>>>\fg$. Since $\fm_\et$ is contained in this fixed 
point set and $\Ad(z)$ acts trivially on it, we see that $t_0^2\z^\et y=y$ for 
all $y\in\fm_\et$.

If $\fm_\et\ne0$, it follows that $t_0^2\z^\et=1$. 
Thus we have $\Ad(\io(t)z)\th(x)=t_0^{2N/\et}\z^N x=(t_0^2\z^\et)^{N/\et}x=x$.
(Here we use that $N\in\et\ZZ$.)
Thus $x$ is contained in the fixed point set of $\Ad(\io(t_0)z)\th:\fg@>>>\fg$,
so that $x\in\fm_N$. We see that for any $N\in\et\ZZ$ we have
$$\fm_N=\fg^{0,2N/\et}_{\un N}.\tag f$$ 
Now (f) also holds when $\fm_\et=0$. (In that case we must have $\io=1$ hence
$\fm=\fm_0$, so that $\fm_N=0$ for $N\ne0$. Moreover, by 3.6(d), $\fm$ is a 
Cartan subalgebra of $\fg_{\un0}$ and $Z=M$. We also have 
$\fg^{0,2N/\et}_{\un N}=0$ for $N\ne0$. Thus, for $N\ne0$, (f) states that $0=0$.
If $N=0$, (f) states that $\fm$ is its own centralizer in $\fg_{\un0}$, which 
is clear.)

From (e),(f) we see that $\fm_N\sub{}^\e\tfl^{\un\vp}_N$ for all $N\in\et\ZZ$ and that 

(g) {\it for $N\in\et\ZZ$ we have $\fm_N={}^\e\tfl^{\un\vp}_N$ if 
and only if the following holds: $\{(\a,n)\in\car_{\un N};n+\la\vp:\a\ra=2N/\et\}$ is equal 
to $\{(0,2N/\et)\}$ if $(0,2N/\et)\in\car_{\un N}$ and is empty if 
$(0,2N/\et)\n\car_{\un N}$;}
\nl

(g${}'$) {\it for $N\in\ZZ-\et\ZZ$ we have $\fm_N={}^\e\tfl^{\un\vp}_N$ (or equivalently
${}^\e\tfl^{\un\vp}_N=0$) if and only if $\{(\a,n)\in\car_{\un N};n+\la\vp:\a\ra=2N/\et\}=\emp$.}

The condition in (g) can be also expressed as follows:

$\{(\a,n)\in\car^*_{\un N};n+\la\vp:\a\ra=2N/\et\}=\emp$ for any $N\in\et\ZZ$.
\nl
The condition in (g${}'$) can be also expressed as follows:

$\{(\a,n)\in\car_{\un N};n+\la\vp:\a\ra=2N/\et\}=\emp$ for any $N\in\ZZ-\et\ZZ$.
\nl
Indeed, it is enough to show that if $N\in\ZZ-\et\ZZ$ and 
$\{(\a,n)\in\car_{\un N};n+\la\vp:\a\ra=2N/\et\}$ then we have automatically $\a\ne0$.
(Assume that $\a=0$. Then $n=2N/\et$ so that $n$ is an odd integer. Since $\fg^0=\fm$
and $\fg^{0,n}_{\un N}\ne0$ we see that ${}^\io_n\fm\ne0$. Using 1.2(d) we deduce that
$n$ is even, a contradiction.)

We see that (b) holds.

From (c) we deduce 
$${}^\e\fp_N^{\un\vp}=\op_{k\in\ZZ;k\ge2fN/\et}({}^\l_k\fg_{\un N})=
\op_{k\in\ZZ,(\a,n)\in\car_{\un N};k\ge2fN/\et,nf+\la\l':\a\ra=k}
(\fg^{\a,n}_{\un N}),$$
hence
$${}^\e\fp^{\un\vp}_N=\op_{(\a,n)\in\car_{\un N};
\la\vp:\a\ra\ge2N/\et-n}(\fg^{\a,n}_{\un N}).\tag h$$
The nilradical ${}^\e\fu_*^{\un\vp}$ of ${}^\e\fp_*^{\un\vp}$ is given by
$${}^\e\fu_N^{\un\vp}=\op_{(\a,n)\in\car^*_{\un N};\la\vp:\a\ra>2N/\et-n}
(\fg^{\a,n}_{\un N}).\tag i$$
We see that for $\vp,\vp'$ in $\EE$ and $N\in\ZZ$, the following two 
conditions are equivalent:

(I)  ${}^\e\fp^{\un\vp}_N={}^\e\fp^{\un\vp'}_N$;

(II) for any $(\a,n)\in\car^*_{\un N}$ we have
$\la\vp:\a\ra+n\ge2N/\et\Lra\la\vp':\a\ra+n\ge2N/\et$.
\nl
For $\vp,\vp'$ in $\EE'$ we say that $\vp\equiv\vp'$ if for any $N\in\ZZ$ and 
any $(\a,n)\in\car^*_{\un N}$ we have
$$(\la\vp:\a\ra+n-2N/\et)(\la\vp':\a\ra+n-2N/\et)>0.$$
This is clearly an equivalence relation on $\EE'$. From the equivalence of 
(I),(II) above, we see that
$$\vp\equiv\vp'\Lra{}^\e\fp^{\un\vp}_N={}^\e\fp^{\un\vp'}_N\qua\frl N\in\ZZ.
\tag j$$
For any $\vp\in\EE'$ we set
$$\align&
I_\vp={}^\e\Ind_{{}^\e\fp^{\un\vp}_\et}^{\fg_\da}(\tC[-\dim\fm_\et])
\in\cq(\fg_\da),\\&
\tI_\vp={}^\e\tInd_{{}^\e\fp^{\un\vp}_\et}^{\fg_\da}(\tC)\in\cq(\fg_\da).\tag k
\endalign$$
Here we regard $\tC$ as an object of 
$\cq(\tfl^{\un\vp}_\et)=\cq(\fm_\et)$, see (b).
Note that in $\ck(\fg_\da)$ we have 
$$\tI_\vp=v^{h(\vp)}I_\vp$$ 
where
$$h(\vp)=\dim{}^\e\fu^{\un\vp}_0+\dim{}^\e\fu^{\un\vp}_\et+\dim\fm_\et=
\dim{}^\e\fu^{\un\vp}_0+\dim{}^\e\fp^{\un\vp}_\et.$$
We show:

(l) {\it If $\vp,\vp'\in\EE'$, $\vp\equiv\vp'$ then $I_\vp=I_{\vp'}$,
$h(\vp)=h(\vp')$.}
\nl
Indeed, in this case we have ${}^\e\fp^{\un\vp}_N={}^\e\fp^{\un\vp'}_N$ for 
all $N\in\ZZ$ (see (j)) and the result follows from the definitions.

\subhead 10.3\endsubhead
We keep the setup of 10.1, 10.2. As in 2.9, for any $N\in\ZZ$ we set
$$\tfl^\ph_N={}_{2N/\et}^\io\fg_{\un N}\text{ if }2N/\et\in\ZZ,\qua \tfl^\ph_N=0\text{ if }2N/\et\n\ZZ.$$
Hence 
$$\tfl^\ph_N=\op_{(\a,n)\in\car_{\un N};n=2N/\et}\fg^{\a,n}_{\un N}.\tag a$$ 
We set $\tfl^\ph=\op_{N\in\ZZ}\tfl^\ph_N$. 

Let 
$$\EE''=\EE-\cup_{N\in\ZZ,(\a,n)\in\car^*_{\un N};n\ne2N/\et}\fH_{\a,n,N}.$$
For $\vp\in\EE$ we show:

(b) {\it We have ${}^\e\tfl^{\un\vp}\sub\tfl^\ph$ if and only if 
$\vp\in\EE''$.}
\nl
Using (a) and 10.2(e) we see that we have $\op_N{}^\e\tfl^{\un\vp}_N\sub\tfl^\ph$ if 
and only if for any $N\in\ZZ$ we have
$$\{(\a,n)\in\car_{\un N};n+\la\vp:\a\ra=2N/\et\}
\sub\{(\a,n)\in\car_{\un N};n=2N/\et\},$$
or equivalently
$$\{(\a,n)\in\car_{\un N};n+\la\vp:\a\ra=2N/\et;n\ne2N/\et\}=\emp,$$
or equivalently
$$\{(\a,n)\in\car^*_{\un N};\la\vp:\a\ra=2N/\et-n\ne0\}=\emp.$$
This is the same as the condition that $\vp\in\EE''$. This proves (b).

\subhead 10.4\endsubhead
We show: 

(a) {\it Let $\fp_*$ be an $\e$-spiral with a splitting $\fm'_{*}$ such that $\fm$ 
is a Levi subalgebra of a parabolic subalgebra of $\fm'=\oplus_{N}\fm'_{N}$ compatible 
with the $\ZZ$-gradings. Then for some $\vp\in\EE$ we have 
$\fp_*={}^\e\fp^{\un\vp}_*$, $\fm'_*={}^\e\tfl^{\un\vp}_*$.}
\nl
We can find $\mu\in Y_{G_{\un0},\QQ}$ such that $\fp_*={}^\e\fp^\mu_*$ 
and $\fm'_{*}={}^\e\tfl^\mu_*$. Let $M'_{0}=e^{\fm'_{0}}$. We choose 
$f\in\ZZ_{>0}$ such that $\l_1:=f\mu\in Y_{G_{\un0}}$. 
We have $\fm'_N={}^{\l_1}_{fN\e}\fg_{\un N}$ for
all $N\in\ZZ$. Let $N\in\et\ZZ$. Since $\fm_{N}\subset\fm'_{N}$, for any $x\in\fm_N$ 
and any $t\in\kk^*$ we have $\Ad(\l_1(t))x=t^{fN\e}x$; we have also 
$\Ad(\io(t))x=t^{2N/\et}x$.

 Hence $\Ad(\l_1(t^2)\io(t^{-f|\et|}))x=x$ for any $x\in\fm_{N}$. 
Since this holds for any $N\in\et\ZZ$, it follows that the image of 
$\l_1^2\io^{-f|\et|}:\kk^*@>>>M'_{0}$ commutes with $M$. Since $M$ 
is a Levi subgroup of a parabolic subgroup of $M'$, the image of $\l_1^2\io^{-f|\et|}$ is 
contained in $\cz_M$, hence in $\cz_M^0=Z$. In particular, the images of $\io$ and 
$\l_1^2\io^{-f|\et|}$ commute with each other, hence the images of $\l_1$ and $\io$ commute 
with each other. It therefore makes sense to write $\l':=2\l_1-f|\et|\io\in Y_{M}$ and we actually 
have $\l'\in Y_Z$. Let $\vp=|\et|\i f\i\l'\in Y_{Z,\QQ}=\EE$. We have
$\vp=|\et|\i(2\mu-|\et|\io)$ hence $\vp+\io=2|\et|\i \mu$ that is, $\mu=\un\vp$
and (a) is proved.

\mpb

Let $\cc'$ be the collection of $\e$-spirals $\fp_*$ such that $\fm_*$ is a 
splitting of $\fp_*$. We show:

(b) {\it $\cc'$ coincides with the collection of $\e$-spirals of the form 
${}^\e\fp_*^{\un\vp}$ with $\vp\in\EE'$.}
\nl
Assume first that $\fp_*\in\cc'$. Using (a) with $\fm'_*=\fm_*$ we see that 
for some $\vp\in\EE$ we have $\fp_*={}^\e\fp_*^{\un\vp}$, 
$\fm_*={}^\e\tfl^{\un\vp}_*$. Using 10.2(b), we see that $\vp\in\EE'$.

Conversely, assume that $\fp_*={}^\e\fp_*^{\un\vp}$ for some $\vp\in\EE'$. 
From 10.2(b) we have $\fm_*={}^\e\tfl_*^{\un\vp}$. Thus $\fp_*\in\cc'$. This 
proves (b).

\mpb

Let $\cc''$ be the collection of $\e$-spirals $\fp_*$ with the following 
property: there exists a splitting $\fm'_*$ of $\fp_*$ such that 
$\fm_N\sub\fm'_N\sub\tfl^\ph_N$ for all $N$. We show:

(c) {\it $\cc''$ coincides with the collection of $\e$-spirals of the form
${}^\e\fp_*^{\un\vp}$ with $\vp\in\EE''$.}
\nl
Assume first that $\fp_*\in\cc''$ and let $\fm'_*$ be a splitting of $\fp_*$
as in the definition of $\cc''$. 
Now $\fm$ is a Levi subalgebra of a parabolic subalgebra of $\fm'=\op_N\fm'_N$ 
compatible with the $\ZZ$-gradings of $\fm$ and $\fm'$. (Indeed, from the 
proof of 3.7(c) we see that there exists $\l\in Y_Z$ such that 
$\fm=\{y\in\tfl^\ph;\Ad(\l(t))y=y\qua\frl t\in\kk^*\}$. Since
$\fm'\sub\tfl^\ph$ we see that 
$\fm=\{y\in\fm';\Ad(\l(t))y=y\qua\frl t\in\kk^*\}$, as required.)
Using (a), we see that for some $\vp\in\EE$
we have $\fp_*={}^\e\fp_*^{\un\vp}$, $\fm'_*={}^\e\tfl_*^{\un\vp}$. Since 
${}^\e\tfl^{\un\vp}_N\sub\tfl^\ph_N$
for all $N$, we see from 10.3(b) that $\vp\in\EE''$.

Conversely, assume that 
$\fp_*={}^\e\fp_*^{\un\vp}$ for some $\vp\in\EE''$. Let 
$m'_*={}^\e\tfl_*^{\un\vp}$. 
From 10.3(b) we see that ${}^\e\tfl^{\un\vp}_N\sub\tfl^\ph_N$ 
for all $N$. We also have $\fm_N\sub{}^\e\tfl^{\un\vp}_N$ for all $N$. Thus 
$\fp_*\in\cc''$. This proves (c).

\mpb

Note that 

(d) {\it if $\fp_*\in\cc''$ then the splitting $\fm'_*$ in the definition
of $\cc''$ is in fact unique.}
\nl
Indeed, assume that $\fm'_*,\tfm'_*$ are 
splittings of $\fp_*$ such that $\fm_N\sub\fm'_N$, $\fm_N\sub\tfm'_N$ for all 
$N$. By 2.7 we can find $u\in U_0$ ($U_0$ as in 2.5) such that
$\Ad(u)\fm'_*=\tfm'_*$. In particular, we have $\Ad(u)\fm'_0=\tfm'_0$. This 
implies $u=1$ since $\fm'_0,\tfm'_0$ are both Levi subalgebras of $\fp_0$ containing 
$\fm_{0}$. Hence $\fm'_*=\tfm'_*$.

\subhead 10.5\endsubhead
Let $\fp_*$, $\fp'_*$ be two $\e$-spirals such that $\fp_N\sub\fp'_N$ for all 
$N$. Let $\fu_*,\fu'_*$ be their nilradicals. Then we have $\fu'_N\sub\fu_N$
for all $N$. In particular, $\fu'_N\sub\fp_N$ for all $N$. We show:

(a) {\it $\op_N\fp_N/\fu'_N$ is a parabolic subalgebra of 
$\fl'=\op_N\fp'_N/\fu'_N$.}
\nl
Let $P_0=e^{\fp_0}$, $P'_0=e^{\fp'_0}$ and let $U_0=U_{P_0},U'_0=U_{P'_0}$.
We have $P_0\sub P'_0$, $U'_0\sub U_0$. Now $\fp_*={}^{\e}\fp^{\mu}_{*}$  and 
$\fp'_*={}^{\e}\fp^{\mu'}_{*}$ for some $\mu=\l/r,\mu'=\l'/r'$, with $\l,\l'\in Y_{G_{\un0}}$ and $r,r'\in\ZZ_{>0}$.
We have $\l(\kk^*)\sub P_0$, $\l'(\kk^*)\sub P'_0$. 
We can find Levi subgroups $\tL_0$, $\tL'_0$ of $P_0,P'_0$ such that 
$\tL_0\sub \tL'_0$. 
By conjugating $\l$ (resp. $\l'$) by an element of $U_0$
(resp. $U'_0$) we can assume that $\l\in\cz_{\tL_0}$, $\l'\in\cz_{\tL'_0}$. Since 
$\cz_{\tL'_0}\sub\cz_{\tL_0}$, we have $\l(t)\l'(t')=\l'(t')\l(t)$ for any $t,t'$ 
in $\kk^*$. 
Hence we have $\fg=\op_{k,k'\in\QQ,i\in\ZZ/m}({}^{\mu,\mu'}_{k,k'}\fg_i)$
where ${}^{\mu,\mu'}_{k,k'}\fg_i={}^\mu_k\fg_i\cap{}^{\mu'}_{k'}\fg_i$. We have 

$\fp_N=\op_{k,k'\in\QQ;k\ge N\e}({}^{\mu,\mu'}_{k,k'}\fg_{\un N})$,

$\fp'_N=\op_{k,k'\in\QQ;k'\ge N\e}({}^{\mu,\mu'}_{k,k'}\fg_{\un N})$,

$\fu'_N=\op_{k,k'\in\QQ;k'> N\e}({}^{\mu,\mu'}_{k,k'}\fg_{\un N})$.
\nl
Since $\fu'_N\sub\fp_N\sub\fp'_N$ we see that

$\fp_N=\op_{k,k'\in\QQ;k\ge N\e,k'\ge N\e}({}^{\mu,\mu'}_{k,k'}\fg_{\un N})$,

$\fu'_N=\op_{k,k'\in\QQ;k\ge N\e,k'>N\e}({}^{\mu,\mu'}_{k,k'}\fg_{\un N})$,
\nl
hence

$\fp_N/\fu'_N\cong \op_{k,k'\in\QQ;k\ge N\e,k'= N\e}({}^{\mu,\mu'}_{k,k'}\fg_{\un N})$.
\nl
This is a subspace of

$\fp'_N/\fu'_N\cong\tfl'_{N}:=\op_{k,k'\in\QQ;k'= N\e}({}^{\mu,\mu'}_{k,k'}\fg_{\un N})$.
\nl
Since $\mu$ and $\mu'$ commute with each other, it makes sense to define $\nu=\mu-\mu'\in Y_{\tL'_{0},\QQ}$.  
Let $\tfl'=\op_N\tfl'_N\sub \fg$. Then $\tL'_0$ acts on $\fl'$ by the $\Ad$-action and $\nu$ induces a $\QQ$-grading 
$\fl'=\op_{k_1\in\QQ}{}^\nu_{k_1}\fl'$. From the definitions we see 
that $\op_N\fp'_N/\fu_N=\op_{k_1\in\QQ;k_1\ge0}({}^\nu_{k_1}\fl')$. Thus (a) 
holds.

\mpb

A similar argument shows:

(b) {\it If $\fm_*$ is a splitting of $\fp_*$, the obvious map $\fm=\op_N\fm_N@>>>
\op_N\fp_N/\fu'_N$ defines an isomorphism of $\fm$ onto a Levi
subalgebra of $\op_N\fp_N/\fu'_N$.}

\subhead 10.6\endsubhead
Let $\vp,\vp'$ in $\EE'$. For any $t\in\QQ$ such that $0\le t\le1$ we set
$$\vp_t:=t\vp+(1-t)\vp'.$$
We assume that there is a unique hyperplane $\fH$ of the form
$\fH=\fH_{\a_0,n_0,N_0}$ for some $(\a_0,n_0,N_0)$ 
with $N_0\in\ZZ$, $(\a_0,n_0)\in\car^*_{\un N_0}$ such that 
$\vp_t\in\fH$ for some $t=s$; this $s$ is necessarily unique 
since $\vp_0\n\fH$. Note however that the triple $(\a_0,n_0,N_0)$ is not uniquely 
determined by $\fH$. We set $\vp''=\vp_s$,
$$\align&\fp_*={}^\e\fp^{\un\vp}_*, \fp'_*={}^\e\fp^{\un\vp'}_*,
\fu_*={}^\e\fu^{\un\vp}_*, \fu'_*={}^\e\fu^{\un\vp'}_*,\\&
\fp''_*={}^\e\fp^{\un\vp''}_*,\tfl''_*={}^\e\tfl^{\un\vp''}_*,
\tfl''=\op_{N\in\ZZ}\tfl''_N,\endalign$$
$$P_0=e^{\fp_0}, P'_0=e^{\fp'_0}, P''_0=e^{\fp''_0}.$$
We show: 

(a) {\it For any $N\in\ZZ$ we have $\fp_N\sub\fp''_N$; hence $P_0\sub P''_0$. 
For any $N\in\ZZ$ we have $\fp'_N\sub\fp''_N$; hence $P'_0\sub P''_0$.}
\nl
Using 10.2(h), we see that to prove the first sentence in (a) it is enough to 
show that for any $(\a,n)\in\car^*_{\un N}$ such that $\la\vp:\a\ra\ge2N/\et-n$
we have also $\la\vp'':\a\ra\ge2N/\et-n$. 
(Since $\vp\in\EE'$, we must have $\la\vp:\a\ra>2N/\et-n$. Since for 
$t\in\QQ$, $s<t\le1$ we have $\vp_t\in\EE'$, it follows that for all such $t$ 
we have $\la\vp_t:\a\ra>2N/\et-n$. Taking the limit as $t\m s$ we obtain
$\la\vp'':\a\ra\ge2N/\et-n$, as required). 
The second sentence in (a) is proved in a similar way.

We show: 

(b) {\it Fix $N\in\ZZ$. If $\fH\ne\fH_{\a,n,N}$ for any $(\a,n)\in\car^*_{\un N}$, 
then we have $\fp_N=\fp''_N=\fp'_N$.}
\nl
We first show that $\fp_N=\fp''_N$. By the equivalence of (I),(II) in 10.2, it 
is enough to show that for any $(\a,n)\in\car^*_{\un N}$ we have
$$\la\vp:\a\ra+n-2N/\et>0\Lra\la\vp'':\a\ra+n-2N/\et>0$$ 
or equivalently that $c_1c_s>0$ where $c_t=\la\vp_t:\a\ra+n-2N/\et$
for $t\in\QQ$ such that $0\le t\le1$. 
From our assumptions we see that $c_t\ne0$ for all $t$.
It follows that either $c_t>0$ for all $t$ or $c_t<0$ for all $t$. In 
particular, $c_1c_s>0$, as required.
The proof of the equality $\fp'_N=\fp''_N$ is entirely similar.

We show:

(c) {\it If $\fH\neq\fH_{\a,n,0}$ for any $(\a,n)\in \car^*_{\un0}$ and 
$\fH\neq\fH_{\a,n,\et}$ for any $(\a,n)\in \car^*_{\da}$ , then $I_\vp=I_{\vp'}$, 
$h(\vp)=h(\vp')$.}
\nl
Indeed, in this case, by (b) we have $\fp_N=\fp'_N$ for $N\in\{0,\et\}$ and the 
equality $I_\vp=I_{\vp'}$ follows from the definitions. 
From $\fp_0=\fp'_0$ we deduce that $\fu_0=\fu'_0$. Using this
and $\fp_\et=\fp'_\et$ we deduce that $h(\vp)=h(\vp')$. This proves (c).

We show:

(d) {\it If $\fH=\fH_{\a_{0},n_{0},0}$ for some $(\a_0,n_0)\in\car^*_{\un0}$ 
but $\fH\neq\fH_{\a,n,\et}$ for any $(\a,n)\in\car^*_\da$, then 
$I_\vp\cong I_{\vp'}$, $h(\vp)=h(\vp')$.}
\nl
By (a) we have $\fp_N\sub\fp''_N,\fp'_N\sub\fp''_N$ for all $N\in\ZZ$ and by (b) 
we have $\fp_\et=\fp''_\et=\fp'_\et$. We can now apply 4.5(b) twice 
to conclude that 
$$I_\vp\cong\op_{j\in J}I_{\vp''}[-2a_j],
I_{\vp'}\cong\op_{j'\in J'}I_{\vp''}[-2a'_{j'}]$$ 
where $a_j$ (resp. $a'_{j'}$) are integers such that
$$\r_{P''_0/P_0!}\bbq=\op_{j\in J}\bbq[-2a_j],
\r_{P''_0/P'_0!}\bbq=\op_{j'\in J'}\bbq[-2a'_{j'}].$$
To show that $I_\vp\cong I_{\vp'}$ it is then enough to show that 
$$\op_{j\in J}\bbq[-2a_j]\cong\op_{j'\in J'}\bbq[-2a'_{j'}]$$ 
or that $\r_{P''_0/P_0!}\bbq\cong\r_{P''_0/P'_0!}\bbq$. This is clear since 
$P''_0/P_0,P''_0/P'_0$ are partial flag manifolds of the reductive quotient 
$L''_{0}$ of $P''_{0}$ with respect to two associate parabolic subgroups. The 
above argument shows also that $\dim U_{P_0}=\dim U_{P'_0}$ hence
$\dim\fu_0=\dim\fu'_0$. Moreover, from (b) we have $\fp_\et=\fp'_\et$; we see 
that $h(\vp)=h(\vp')$. This proves (d).

\mpb

For $N\in\ZZ$ let $\fq_N$ (resp. $\fq'_N$) be the image of $\fp_N$ (resp.
$\fp'_N$) under the obvious projection $\fp''_N@>>>\tfl''_N$. From 10.5(a),(b),
we see that $\fq=\op_N\fq_N$,  $\fq'=\op_N\fq'_N$  are parabolic subalgebras of
$\tfl''$ and $\fm$ is a Levi subalgebra of both $\fq$ and $\fq'$. We show:

(e) {\it If $\fH=\fH_{\a_0,2,\et}$ for some $(\a_0,2)\in\car^*_\da$, then 
$I_\vp\cong I_{\vp'}$ and $h(\vp)=h(\vp')$.}
\nl
By (a) we have $\fp_N\sub\fp''_N,\fp'_N\sub\fp''_N$ for all $N\in\ZZ$. Since 
$\fH_{\a_0,2,\et}$ is the unique hyperplane (as in 10.1)
on which $\vp''$ lies, we have 
$\vp''\in\EE''$. Hence by 10.3(b) we have $\tfl''_N\sub\tfl^\ph_N$ for all 
$N\in\ZZ$. In particular, the $\ZZ$-grading of $\tfl''$ is $\et$-rigid and
$e\in\ovsc\tfl''_\et$, so that $\ovsc\fm_\et\sub\ovsc\tfl''_\et$.
Let $A\in\cq(\tfl''_\et)$ be the simple perverse sheaf on $\tfl''_\et$ such that
the support of $A$ is $\tfl''_\et$ and $A|_{\ovsc\fm_\et}$ is
equal up to shift to $\tC|_{\ovsc\fm_\et}$.
Applying (twice) 1.8(b) and the transitivity formula 4.2(a) we deduce
$${}^\e\Ind_{\fp_\et}^{\fg_\da}(\tC)=
{}^\e\Ind_{\fp''_\et}^{\fg_\da}(\ind_{\fq_\et}^{\tfl''_\et}(\tC))\cong
\op_j{}^\e\Ind_{\fp''_\et}^{\fg_\da}(A)[-2s_j][\dim\fm_\et-\dim\tfl''_\et],$$
$${}^\e\Ind_{\fp'_\et}^{\fg_\da}(\tC)=
{}^\e\Ind_{\fp''_\et}^{\fg_\da}(\ind_{\fq'_\et}^{\tfl''_\et}(\tC))\cong
\op_j{}^\e\Ind_{\fp''_\et}^{\fg_\da}(A)[-2s_j][\dim\fm_\et-\dim\tfl''_\et],$$
where $(s_j)$ is a certain finite collection of integers. It follows that
$${}^\e\Ind_{\fp_\et}^{\fg_\da}(\tC)\cong{}^\e\Ind_{\fp'_\et}^{\fg_\da}(\tC).$$
This shows that $I_\vp,I_{\vp'}$ are isomorphic. It remains to show that 
$h(\vp)=h(\vp')$, for which it suffices to show the following two equalities:
$$\dim\fu_0=\dim\fu'_0,\tag f$$ 
$$\dim\fu_\et=\dim\fu'_\et.\tag g$$
Since $\fp_0$ and $\fp'_0$ are associate parabolic subalgebras of $\fp''_0$ 
sharing the same Levi subalgebra $\fm_0$, their unipotent radicals $\fu_0$ and
$\fu'_0$ have the same dimension. This proves (f). 

Now we show (g). By 10.2(i), we have
$$\dim\fu_\et-\dim\fu'_\et=\sum_{\a\in X_Z;\la\vp:\a\ra>0,\la\vp':\a\ra<0}
\dim\fg^{\a,2}_\da-\sum_{\a\in X_{Z};\la\vp:\a\ra<0,\la\vp':\a\ra>0}
\dim\fg^{\a,2}_\da.$$
We want to show the above difference is zero. For this, it suffices to show 
that $\dim\fg^{\a,2}_\da=\dim\fg^{-\a,2}_\da$ for any $\a\in X_Z$. Note that 
$\fg^{\a,2}_\da=\tfl^\phi_\et\cap\fg^\a$,  which is the $\a$-weight space of 
$Z$ on $\tfl^\phi_\et$. By the property of the $\fs\fl_2$-action on 
$\tfl^\phi$ given by the triple $\phi=(e,h,f)$, we have that $\ad(e)$ induces 
an isomorphism $\tfl^\phi_{-\et}\cong\tfl^\phi_\et$. Since $Z$ commutes with 
$e$, this isomorphism restricts to an isomorphism of $\a$-weight spaces 
$\tfl^\phi_{-\et}\cap\fg^\a\cong\tfl^\phi_\et\cap\fg^\a$, that is 
$\fg^{\a,-2}_{-\da}\cong\fg^{\a,2}_\da$. Therefore 
$\dim\fg^{\a,2}_\da=\dim\fg^{\a,-2}_{-\da}=\dim\fg^{-\a,2}_\da$, as desired. This
proves (g) and finishes the proof of (e).

\subhead 10.7\endsubhead
Let
$$\align&\ovsc\EE=\EE-\cup_{(\a,n)\in\car^*_\d;n\ne2}\fH_{\a,n,\et}\\&=
\{\vp\in\EE;\la\vp:\a\ra+n-2\ne0\qua\frl(\a,n)\in\car^*_\da
\text{ such that }n\ne2\}.\endalign$$
Note that $\EE'\sub\EE''\sub\ovsc\EE$.
For ${}'\vp,{}''\vp$ in $\ovsc\EE$ we say that ${}'\vp\si{}''\vp$ if for any
$(\a,n)\in\car^*_\da$ such that $n\ne2$ we have
$$(\la{}'\vp:\a\ra+n-2)(\la{}''\vp:\a\ra+n-2)>0.$$
This is an equivalence relation on $\ovsc\EE$. We show:

(a) {\it Assume that ${}'\vp,{}''\vp$ in $\EE'$ satisfy ${}'\vp\si{}''\vp$. 
Then $I_{{}'\vp}\cong I_{{}''\vp}$ and $h({}'\vp)=h({}''\vp)$; hence
$\tI_{{}'\vp}\cong\tI_{{}''\vp}$.}
\nl
By a known property of hyperplane arrangements, we can find a sequence 
\lb
$\vp_0,\vp_1,\do,\vp_k$ in $\EE'$ such that $\vp_0={}'\vp,\vp_k={}''\vp$ and 
such that for any $j\in\{0,1,\do,k-1\}$ one of the following holds:

$\vp_j\equiv\vp_{j+1}$;

$\vp_j=\vp,\vp_{j+1}=\vp'$ are as in 10.6(c);

$\vp_j=\vp,\vp_{j+1}=\vp'$ are as in 10.6(d);

$\vp_j=\vp,\vp_{j+1}=\vp'$ are as in 10.6(e).
\nl
Using 10.2(l) or 10.6(c) or 10.6(d) or 10.6(e) we see that for 
$j\in\{0,1,\do,k-1\}$ we have $I_{\vp_j}\cong I_{\vp_{j+1}}$ and 
$h(\vp_j)=h(\vp_{j+1})$. This proves (a).

\mpb

The set of equivalence classes on $\ovsc\EE$ for $\si$ is denoted by 
$\un{\ovsc\EE}$; it is a finite set. 

\subhead 10.8\endsubhead
For $\vp,\vp'$ in $\EE$ we set  
$$\align&\t(\vp,\vp')=
\sum_{(\a,n)\in\car^*_\da;(\la\vp:\a\ra+n-2)(\la\vp':\a\ra+n-2)<0}
\dim\fg^{\a,n}_\da\\&
-\sum_{(\a,n)\in\car^*_{\un0};(\la\vp:\a\ra+n)(\la\vp':\a\ra+n)<0}
\dim\fg^{\a,n}_{\un0}.\tag a\endalign$$
Using 10.2(i) we see that when $\vp,\vp'$ are in $\EE'$, then
$$\t(\vp,\vp')=-\dim\fra{{}^\e\fu_0^{\un\vp}+{}^\e\fu_0^{\un\vp'}}
{{}^\e\fu_0^{\un\vp}\cap{}^\e\fu_0^{\un\vp'}}
+\dim\fra{{}^\e\fu_\et^{\un\vp}+{}^\e\fu_\et^{\un\vp'}}
{{}^\e\fu_\et^{\un\vp}\cap{}^\e\fu_\et^{\un\vp'}}.\tag b$$

\subhead 10.9\endsubhead
We define $G_\ph,M_\ph$ as in 3.6. We show that 

(a) {\it the obvious map
$M_\ph/M_\ph^0@>>>(G_{\un0}\cap G_\ph)/(G_{\un0}\cap G_\ph)^0$
is an isomorphism.}
\nl
Recall that $\ph=(e,h,f)$ with $e\in\ovsc\fm_\et$. Let $U$ (resp. $U'$) be the 
unipotent radical of $G(e)^0$ (resp. $(M\cap G(e))^0$). We have 
$G(e)=G_\ph U$ (semidirect product). Taking fixed point sets of $\vt$ we obtain
$G_{\un0}\cap G(e)=(G_{\un0}\cap G_\ph)(G_{\un0}\cap U)$ (semidirect product).
We have $M\cap G(e)=M_\ph U'$ (semidirect product).
Taking fixed point set of $\io(t)$ for all $t\in\kk^*$ we obtain
$M_0\cap G(e)=(M_0\cap M_\ph)(M_0\cap U')=M_\ph(M_0\cap U')$ (semidirect 
product). (We have used that $M_\ph\sub M_0$.)  
It follows that we have canonically
$$(M_0\cap G(e))/(M_0\cap G(e))^0=M_\ph/M_\ph^0,$$
$$(G_{\un0}\cap G(e))/(G_{\un0}\cap G(e))^0
=(G_{\un0}\cap G_\ph)/(G_{\un0}\cap G_\ph)^0.$$
It remains to use that
$$(M_0\cap G(e))/(M_0\cap G(e))^0=(G_{\un0}\cap G(e))/(G_{\un0}\cap G(e))^0,$$ 
see 3.8(a).

\mpb

Recall from 3.6 that $Z$ is a maximal torus of $(G_\ph\cap G_{\un0})^0$. Let 
$H$ be the normalizer of $Z$ in $(G_\ph\cap G_{\un0})^0$. Let 
$H'=\{g\in G_{\un0};\Ad(g)M=M,\Ad(g)\fm_k=\fm_k\qua\frl k\in\ZZ\}$. Note that 
$M_0$ is a normal subgroup of $H'$. Hence the groups $H/Z$, $H'/M_0$ are 
defined. We show:

(b) $H\sub H'$.
\nl
Let $g\in(G_\ph\cap G_{\un0})^0$ be such that $gZg\i=Z$. Let 
$$(M',M'_0,\fm',\fm'_*)=(gMg\i,gM_0g\i,\Ad(g)\fm,\Ad(g)\fm_*).$$
Note that $\cz_{M'}^0=\cz_M^0=Z$. Repeating the argument in the last paragraph
in the proof of 3.6(a), we see that
$$(M',M'_0,\fm',\fm'_*)=(M,M_0,\fm,\fm_*).$$ 
(The argument is applicable since $g\in G_\ph\cap G_{\un0}$.) We see that 
$g\in H'$; this proves (b).

We show:

(c) $H\cap M_0=Z$.
\nl
From the injectivity of the map in (a) we see that
$M_\ph\cap(G_{\un0}\cap G_\ph)^0=M_\ph^0$. Hence we have 
$$M_0\cap(G_{\un0}\cap G_\ph)^0\sub(M_0\cap G_\ph)\cap(G_{\un0}\cap G_\ph)^0
\sub M_\ph\cap(G_{\un0}\cap G_\ph)^0=M_\ph^0,$$
so that 
$$M_0\cap(G_{\un0}\cap G_\ph)^0\sub M_\ph^0.$$
The opposite inclusion is also true since $M_\ph\sub M_0$ and
$M_\ph^0\sub G_{\un0}\cap G_\ph$. It follows that
$$M_0\cap(G_{\un0}\cap G_\ph)^0=M_\ph^0=Z.$$ 
The last equality is because $e$ is distinguished in $\fm$. 
Now $H\cap M_0$ is the normalizer of $Z$
in $M_0\cap(G_\ph\cap G_{\un0})^0$ that is, the normalizer of $Z$ in $Z$. We 
see that $H\cap M_0=Z$. This proves (c).

We show:

(d) $H'=M_0H$.
\nl
Since $H\sub H'$, $M_0\sub H'$, we have $M_0H\sub H'$. Now let $g\in H'$. We 
show that $g\in M_0H$. Let $\ph'=(\Ad(g)e,\Ad(g)h,\Ad(g)f)$. We have
$$\Ad(g)e\in\ovsc\fm_\et,\Ad(g)h\in\fm_0,\Ad(g)f\in\fm_{-\et}.$$ 
Since both $\Ad(g)e,e$ are in $\ovsc\fm_\et$, we can find $g_1\in M_0$ such that
$\Ad(g')\Ad(g)e=e$. Replacing  $g$ by $g_1g$ we can assume that we have
$\Ad(g)e=e$. Using \cite{\GRA, 3.3} for $J^M$, we can find $g_2\in M_0$ such 
that 
$$(\Ad(g_2)\Ad(g)e,\Ad(g_2)\Ad(g)h,\Ad(g_2)\Ad(g)f)=(e,h,f).$$ 
We have $g_2g\in G_\ph$. Replacing $g$ by $g_2g$ we can assume that 
$g\in G_{\un0}\cap G_\ph$. Using the surjectivity of the map in (a) we see that
$$G_{\un0}\cap G_\ph\sub M_\ph(G_{\un0}\cap G_\ph)^0.$$ 
Thus we can write $g$ in the form $g_3g'$ with $g_3\in M_\ph$, 
$g'\in(G_{\un0}\cap G_\ph)^0$. Replacing $g$ by $g_3\i g$ we see that we can 
assume that $g\in(G_{\un0}\cap G_\ph)^0$. Since $\Ad(g)M=M$, we see that
$\Ad(g)Z=Z$. Thus $g\in H$. This proves (d).

\mpb

From (b),(c),(d) we see that 

(e) {\it the inclusion $H\sub H'$ induces an isomorphism $H/Z@>\si>>H'/M_0$.
In particular, $M_0$ is the identity component of $H'$.}

\subhead 10.10\endsubhead
Let $g\in H'$ (notation of 10.9). Then $\Ad(g)$ restricts to an isomorphism
$\fm_\et@>\si>>\fm_\et$. Let $\tC'=\Ad(g)^*\tC$, a simple perverse sheaf in
$\cq(\fm_\et)$. We show:

(a) $\tC'\cong\tC$.
\nl
Using 10.9(e) we can assume that $g\in H$ (notation of 10.9). Since $\tC,\tC'$ 
are intersection cohomology complexes attached to $M_0$-equivariant 
irreducible local systems on $\fm_\et$, they correspond to irreducible 
representations of 
$$(M_0\cap G(e))/(M_0\cap G(e))^0=M_\ph/M_\ph^0.$$
Hence it is enough to show that $\Ad(g)$ induces the identity automorphism of 
$M_\ph/M_\ph^0$. Using 10.9(a), we see that it is enough to show that $\Ad(g)$ 
induces the identity automorphism of 
$(G_{\un0}\cap G_\ph)/(G_{\un0}\cap G_\ph)^0$. This is obvious 
since \lb 
$g\in(G_{\un0}\cap G_\ph)^0$. This proves (a).

\subhead 10.11\endsubhead   
Let $\vp,\vp'\in\EE'$. Recall that $P_0=e^{{}^\e\fp_0^{\un\vp}}\sub G_{\un0}$, 
$P'_0=e^{{}^\e\fp_0^{\un\vp'}}\sub G_{\un0}$ are parabolic subgroups of 
$G_{\un0}$ with a common Levi subgroup $M_0$. Let $U_0=U_{P_0}$, 
$U'_0=U_{P'_0}$.
Let $X$ be the set of all $g\in G_{\un0}$ such that $\Ad(g)\fp_*^{\un\vp}$ and 
$\fp_*^{\un\vp'}$ have a common splitting (as in 6.3). 
Note that $X$ is a union of $(P'_0,P_0)$-double cosets in $G_{\un0}$. We show:

(a) {\it We have $H'\sub X$ (notation of 10.9). Let 
$j:H'/M_0@>>>P'_0\bsl X/P_0$ be the map induced by the inclusion $H'@>>>X$.
Then $j$ is a bijection.}
\nl
If $g\in H'$ then $g\in G_{\un0}$ and $\Ad(g)\fm_*=\fm_*$ hence $\fm_*$ is a 
common splitting of $\Ad(g)\fp_*^{\un\vp}$ and $\fp_*^{\un\vp'}$. Thus we have 
$g\in X$, so that the inclusion $H'\sub X$ holds.
Let $g\in X$. Then $g\in G_{\un0}$ and $\Ad(g)\fp_*^{\un\vp}$, 
$\fp_*^{\un\vp'}$ have a common splitting $\fm'_*$. Then $\Ad(g)\fm_*,\fm'_*$ 
are splittings of $\Ad(g)\fp_*^\vp$ hence, by 2.7(a), we have 
$\Ad(gug\i)\Ad(g)\fm_*=\fm'_*$ for some $u\in U_0$. Moreover, $\fm_*,\fm'_*$ 
are splittings of $\fp_*^{\un\vp'}$ hence, by 2.7(a), we have 
$\Ad(u')\fm_*=\fm'_*$ for some $u'\in U'_0$. It follows that 
$\Ad(gu)\fm_*=\Ad(u')\fm_*$ hence 
$u'{}\i gu\in H'$. Since $U_0\sub P_0$, $U'_0\sub P'_0$, we see that $j$ is 
surjective. It remains to show that $j$ is injective. Let $g,g'$ be elements 
of $H'$ such that $g'=p'_0gp_0$ for some $p_0\in P_0,p'_0\in P'_0$. We must 
only show that $g'\in gM_0$. Let $NM_0$ be the normalizer of $M_0$ in 
$G_{\un0}$. It is enough to show that the obvious map 
$NM_0/M_0@>>>P'_0\bsl G_{\un0}/P_0$ is injective. 
This is a well known property of parabolic subgroups and their Levi subgroups 
in a connected reductive group. This completes the proof of (a).

\mpb

Let $\cw=H/Z$ (notation of 10.9). Let $w\in\cw$ and let $g\in H$ be a 
representative of $w$. Now $\Ad(g)$ restricts to an automorphism of $Z$ which 
depends only on $w$; this induces an isomorphism $Y_Z@>\si>>Y_Z$ and, by 
extension of scalars, a vector space isomorphism $\EE@>\si>>\EE$ denoted by 
$\vp_1\m w\vp_1$. For any $(\a,n,i)\in X_Z\T\ZZ\T\ZZ/m$, $\Ad(g)$ defines an 
isomorphism $\fg^{\a,n}_i@>\si>>\fg^{{}^w\a,n}_i$ where ${}^w\a\in X_Z$ is 
given by ${}^w\a(z)=\a(w\i(z))$; hence for any $i$, $(\a,n)\m({}^w\a,n)$ is a 
bijection $\car_i@>\si>>\car_i$. Moreover, for any $N\in\ZZ$ and any 
$(\a,n)\in\car^*_{\un N}$, $w:\EE@>>>\EE$ restricts to a bijection from the 
affine hyperplane $\fH_{\a,n,N}$ to the affine hyperplane $\fH_{{}^w\a,n,N}$.
It follows that $w:\EE@>>>\EE$ restricts to a bijection $\EE'@>\si>>\EE'$.

We show: 

(b) {\it For any $\vp\in\EE'$ we have 
$\Ad(g)({}^\e\fp_*^{\un\vp})={}^\e\fp_*^{\un{w\vp}}$.}   
\nl
From 10.2(h) we have  
$$\Ad(g)({}^\e\fp_N^{\un\vp})=\op_{(\a,n)\in\car_{\un N};\la\vp:\a\ra\ge2N/\et-n}
\Ad(g)\fg^{\a,n}_{\un N}
=\op_{(\a,n)\in\car_{\un N};\la\vp:\a\ra\ge2N/\et-n}(\fg^{{}^w\a,n}_{\un N}),$$
$${}^\e\fp_N^{\un{w\vp}}=\op_{(\a,n)\in\car_{\un N};\la w\vp:\a\ra\ge2N/\et-n}
(\fg^{\a,n}_{\un N})
=\op_{(\a',n)\in\car_{\un N};\la w\vp:{}^w\a'\ra\ge2N/\et-n}
(\fg^{{}^w\a',n}_{\un N}).$$ 
It remains to use that $\la\vp:\a\ra=\la w\vp:{}^w\a\ra$.

\subhead 10.12\endsubhead
For $\vp_1,\vp_2$ in $\EE$ we set
$$[\vp_1|\vp_2]=(1-v^2)^{-\dim Z}\sum_{w\in\cw}v^{\t(\vp_2,w\vp_1)}\in\QQ(v).
\tag a$$
(Here $v$ is an indeterminate and $\t(\vp_2,w\vp_1)\in\ZZ$ is as in 10.8.) 
When $\vp_1,\vp_2$ are in $\EE'$ we have:
$$\sum_{j\in\ZZ}d_j(\fg_\da;\tI_{\vp_1},D(\tI_{\vp_2}))v^{-j}=[\vp_1|\vp_2].
\tag b$$
This can be deduced from 6.4 as follows. The set $X$ in 6.4 is described in 
our case in 10.11(a) in terms of the group $H'/M_0$ which, in turn, is 
identified in 10.9(e) with $\cw=H/Z$; the integers $\t(g)$ in 6.4 are 
identified with the integers $\t(\vp_2,w\vp_1)\in\ZZ$ by 10.11(b). Finally, the 
set $X'$ in 6.4 coincides with $X$ in 6.4 by 10.10(a).

From the definitions we have $\t(\vp_2,w\vp_1)=\t(w\i\vp_2,\vp_1)=
\t(\vp_1,w\i\vp_2)$ for any $w\in\cw$. It follows that 
$$[\vp_1|\vp_2]=[\vp_2|\vp_1]. \tag c$$

\mpb

Let $\boc_1$ (resp. $\boc_2$) be the equivalence class for $\si$ in $\ovsc\EE$
that contains $\vp_1$ (resp. $\vp_2$). Using 10.7(a) we see that the right hand
side of (b) depends only on $\boc_1$, $\boc_2$ and not on the specific elements
$\vp_1\in\boc_1\cap\EE'$, $\vp_2\in\boc_2\cap\EE'$. Hence for $\boc_1,\boc_2$ 
in $\un{\ovsc\EE}$ we can set $[\boc_1|\boc_2]=[\vp_1|\vp_2]\in\QQ(v)$ 
for any $\vp_1\in\boc_1\cap\EE'$, $\vp_2\in\boc_2\cap\EE'$.

\subhead 10.13\endsubhead
Let $\fB$ (resp. ${}^\x\fB$) be the set of (isomorphism classes of) simple
perverse sheaves in $\cq(\fg_\da)$ (resp. ${}^\x\cq(\fg_\da)$).
For any $G_{\un0}$-orbit $\co$ in $\fg_\da^{nil}$ let $\fB_\co$ be the set of 
all $B\in\fB$ such that the support of $B$ is equal to the closure of $\co$.
 We have $\fB=\sqc_\co\fB_\co$. We define a map $\k:\fB@>>>\NN$ by
$\k(B)=\dim\co$ where $B\in\fB_\co$.

\subhead 10.14\endsubhead
Let $\VV'$ be the $\QQ(v)$-vector space with basis
$\{\tT_\boc;\boc\in\un{\ovsc\EE}\}$. On $\VV'$ we have a unique pairing
$(:):\VV'\T\VV'@>>>\QQ(v)$ which is $\QQ(v)$-linear in the first argument,
$\QQ(v)$-antilinear in the second argument (for $f\m\baf$) and such that for
$\boc_1,\boc_2$ in $\un{\ovsc\EE}$ we have 
$(\tT_{\boc_1}:\tT_{\boc_2})=[\boc_1|\boc_2]$ (see 10.12). 

Setting

$\fR_l=\{x\in\VV';(x:x')=0\qua\frl x'\in\VV'\}$,

$\fR_r=\{x\in\VV';(x':x)=0\qua\frl x'\in\VV'\},$
\nl
we state:

\proclaim{Lemma 10.15} We have $\fR_l=\fR_r$.
\endproclaim
The proof is given in 10.17.

\subhead 10.16\endsubhead
We define a $\QQ$-linear map 
$$\ti\g:\VV'@>>>\QQ(v)\ot_\ca{}^\x\ck(\fg_\da)$$ 
by $\tT_\boc\m\tI_\vp$ where 
$\vp$ is an element of $\boc\cap\EE'$. Now $\ti\g$ is well-defined by 10.7(a) 
and is surjective by 8.4(b), 10.4(b). We define a pairing
$$(\QQ(v)\ot_\ca\ck(\fg_\da))\T(\QQ(v)\ot_\ca\ck(\fg_\da))@>>>\QQ((v))$$
(denoted by $(:)$) by requiring that it is $\QQ(v)$-linear in the first 
argument, $\QQ(v)$-antilinear in the second argument (for $f\m\baf$) and that
its restriction 
$$\ck(\fg_\da)\T\ck(\fg_\da)@>>>\QQ((v))$$
is the same as the restriction of the pairing in 4.4(c).
This restricts to a pairing
$$(\QQ(v)\ot_\ca{}^\x\ck(\fg_\da))\T(\QQ(v)\ot_\ca{}^\x\ck(\fg_\da))@>>>\QQ((v))
$$
(denoted again by $(:)$). We show that for $b,b'$ in $\VV'$ we have
$$(\ti\g(b):\ti\g(b'))=(b:b')\tag a$$
We can assume that $b=\tT_\boc,b'=\tT_{\boc'}$ with $\boc,\boc'$ in 
$\un{\ovsc\EE}$. We must show that $\{\tI_\vp,D(\tI_{\vp'})\}=[\vp|\vp']$
where $\vp\in\boc\cap\EE',\vp'\in\boc'\cap\EE'$. This follows from 10.12(b).
This proves (a).

\subhead 10.17\endsubhead
Let ${}'\fR_l=\{z\in\QQ(v)\ot_\ca{}^\x\ck(\fg_\da);(z:z')=0\qua\frl z'\in
\QQ(v)\ot_\ca{}^\x\ck(\fg_\da)\}$,

${}'\fR_r=\{z\in\QQ(v)\ot_\ca{}^\x\ck(\fg_\da);(z':z)=0\qua\frl z'\in
\QQ(v)\ot_\ca{}^\x\ck(\fg_\da)\}$.
\nl
We show:

(a) ${}'\fR_l={}'\fR_r=0$.
\nl
Now $\QQ(v)\ot_\ca{}^\x\ck(\fg_\da)$ has a $\QQ(v)$-basis formed by 
${}^\x\fB=\{B_1,B_2,\do,B_r\}$. From 
0.12 we see that $(B_j:B_{j'})\in\d_{j,j'}+v\NN[[v]]$ for $j,j'$ in $[1,r]$. 

Assume that $\b=\sum_{j\in[1,r]}f_jB_j\in{}'\fR_l$ where 
$f_j\in\QQ(v)$ are not all zero. We must show that this is a contradiction. We
can assume that $f_j\in\QQ[v]$ for all $j$ and $f_{j_0}-c_0\in v\QQ[v]$ for
some $j_0\in[1,r]$ and some $c_0\in\QQ-\{0\}$. Then 
$0=(\b:B_{j_0})\in c_0+v\QQ[[v]]$, a contradiction. This proves that
${}'\fR_l=0$.

Next we assume that $\b=\sum_{j\in[1,r]}f_jB_j\in{}'\fR_r$ where 
$f_j\in\QQ(v)$ are not all zero. We must show that this is a contradiction. We
can assume that $\baf_j\in\QQ[v]$ for all $j$ and 
$\baf_{j_0}-c_0\in v\QQ[v]$ for some $j_0\in[1,r]$ and some 
$c_0\in\QQ-\{0\}$. Then 
$0=(B_{j_0}:\b)\in c_0+v\QQ[[v]]$, a contradiction. This proves that
${}'\fR_r=0$. This proves (a).

We show:

(b) $\fR_l=\ti\g\i({}'\fR_l)$.
\nl
Let $x\in\fR_l$. From 10.16(a) we see that $(\ti\g(x):\ti\g(x'))=0$ for any
$x'\in\VV'$. Since $\ti\g$ is surjective, it follows that
$(\ti\g(x):z')=0$ for any $z'\in \QQ(v)\ot_\ca{}^\x\ck(\fg_\da)$. Thus,
$\ti\g(x)\in{}'\fR_l$.
Conversely, assume that $x\in\VV'$ and $\ti\g(x)\in{}'\fR_l$.
From 10.16(a) we see that for any $x'\in\VV'$ we have
$(x:x')=(\ti\g(x):\ti\g(x'))=0$. Thus $x\in\fR_l$. This proves (b).

An entirely similar proof shows that

(c) $\fR_r=\ti\g\i({}'\fR_r)$.
\nl
From (a),(b),(c) we see that $\fR_l=\fR_r=\ti\g\i(0)$. This proves Lemma 10.15.

\subhead 10.18. Definition of $\VV$\endsubhead
We define $\VV=\VV'/\fR_l=\VV'/\fR_r$ (see 10.15). Note that $(:)$ on $\VV'$ 
induces a pairing $\VV\T\VV@>>>\QQ(v)$ (denoted again by $(:)$) which is
$\QQ(v)$-linear in the first argument, $\QQ(v)$-antilinear in the second 
argument (for $f\m\baf$).

\subhead 10.19\endsubhead
From the proof of 10.15 we see that
$\ti\g$ induces a $\QQ(v)$-linear isomorphism 
$$\g:\VV@>\si>>\QQ(v)\ot_\ca{}^\x\ck(\fg_\da)\tag a$$ 
and that for $b,b'$ in $\VV$ we have
$$(\g(b):\g(b'))=(b:b').\tag b$$
From (a) we deduce the following result:

\proclaim{Proposition 10.20} The number of simple perverse sheaves (up to 
isomorphism) in ${}^\x\cq(\fg_\da)$ is equal to $\dim_{\QQ(v)}\VV$.
\endproclaim 

\subhead 10.21\endsubhead
We define a $\QQ$-linear involution $\bar{}:\VV'@>>>\VV'$ by 
$\ov{f\tT_\boc}=\baf\tT_\boc$ for any $f\in\QQ(v)$, 
$\boc\in\un{\ovsc\EE}$; here $\baf$ is as in 0.12. We show:

(a) {\it For any $x,x'$ in $\VV'$ we have $(x:x')=(\bar x':\bar x)$.}
\nl
We can assume that $x=\tT_\boc,x'=\tT_{\boc'}$ for some $\boc,\boc'$ in
$\un{\ovsc\EE}$. We must show that 
$$[\boc|\boc']=[\boc'|\boc],$$
which follows directly from 10.12(c). This proves (a).  

\mpb

We show:

(b) {\it $\bar{}:\VV'@>>>\VV'$ preserves $\fR_l=\fR_r$ (see 10.15) hence it 
induces an $\QQ(v)$-semilinear involution $\bar{}:\VV@>>>\VV$ (with respect to 
$f\m\baf$).}
\nl
Assume that $x\in\fR_l$. We have $(x:x')=0$ for all $x'\in\VV$. Hence, by (a),
we have $(\bar x':\bar x)=0$ for all $x'\in\VV$. Since $\bar{}:\VV'@>>>\VV'$
is surjective, it follows that $\bar x\in\fR_r$. This proves (b).

\subhead 10.22\endsubhead
Now let $\et_1\in\ZZ-\{0\}$ be such that $\un\et_1=\da$. 
Recall that in 10.1 we have fixed $\x\in\un\fT_\et$ and a representative 
$\dot\x=(M,M_0,\fm,\fm_*,\tC)\in\fT_\et$ for $\x$. 
Let $\dot\x_1=(M,M_0,\fm,\fm_{(*)},\tC)\in\fT_{\et_1}$ be as in 3.9 and let
$\x_1\in\un\fT_{\et_1}$ be the $G_{\un0}$-orbit of $\dot x_1$.
Note that the $\QQ$-vector space $\EE$ defined as in 10.1 in terms of 
$\dot\x$ is the same as $\EE$ defined in terms of $\dot\x_1$. 
Moreover the subset $\ovsc\EE$, the equivalence relation $\si$ on it, and the set
$\un{\ovsc\EE}$ defined as in 10.7 in terms of $\dot\x$ are the same as the analogous objects 
defined in terms of $\dot\x_1$. Also, the pairings $\t:\EE\T\EE@>>>\ZZ$ and
$[?|?]:\EE\T\EE@>>>\QQ(v)$ defined in 10.8 and 10.12 in terms of $\dot\x$ are the same as
those defined in terms of $\dot\x_1$.

The subset $\EE'$ of $\EE$ defined as in 10.1 in terms of 
$\dot\x$ is not in general the same as the analogous subset $\EE'_1$ of $\EE$
defined in terms of $\dot\x_1$. 
However, for $\boc_1,\boc_2$ in
$\un{\ovsc\EE}$, the quantity $[\boc_1|\boc_2]in\QQ(v)$ defined in 10.12 in terms of
$\dot\x$ is the same as that defined in terms of $\dot\x_1$. (It is equal to
$[\vp_1|\vp_2]$ for any $\vp_1\in\boc_1\cap\EE'\cap\EE'_1$, $\vp_2\in\boc_2\cap\EE'\cap\EE'_1$.
Hence the vector space $\VV'$, the pairing $(:)$ on it and its quotient $\VV$ defined in 10.14
and 10.18 in terms of $\dot\x$ are the same as those defined in terms of $\dot\x_1$. 
The involutions $\bar{}:\VV'@>>>\VV'$, $\bar{}:\VV@>>>\VV$  defined in 10.21 
in terms of $\dot\x$ are the same as those defined in terms of $\dot\x_1$. 

\head 11. The $\ca$-lattice $\VV_\ca$\endhead
In this section we give a combinatorial definition of an $\ca$-lattice 
$\VV_\ca$ in the $\QQ(v)$-vector space $\VV$
and a signed basis $\BB'$ of it. It 
turns out that the $\ZZ/2$-orbits on $\BB'$ for the $\ZZ/2$-action $b\m-b$
are in natural bijection with the simple perverse sheaves in the block
${}^\x\cq(\fg_\da)$.

\subhead 11.1\endsubhead
In this section (except in 11.5, 11.6, 11.13) we preserve the setup and notation of 
10.1, 10.2 and assume that $\vp\in\EE''$ (see 10.3). Let
$$\fg^\ph=\{y\in\fg;[y,e]=0,[y,h]=0,[y,f]=0\}.$$
Let $\fz$ be the centre of $\fm$. Note that $\fm\cap\fg^\ph=\fz$
(since $e$ is distinguished in $\fm$) and $\fz$ is a Cartan subalgebra of the 
reductive Lie algebra $\fg^\ph$. (This has already been proved at level of 
groups in 3.6.) For any $\a\in X_Z$ let
$$(\fg^\ph_{\un0})^\a=\fg^\a_{\un0}\cap\fg^\ph.$$
Let 
$$\fg^\ph_{\un0,\vp}=\op_{\a\in X_Z;\la\vp:\a\ra=0}(\fg^\ph_{\un0})^\a.$$
This is a Levi subalgebra (containg $\fz$) of a parabolic subalgebra of 
$\fg^\ph_{\un0}$.
Let $\cb$ be the variety of Borel subalgebras of $\fg^\ph_{\un0,\vp}$, let 
$d(\vp)=\dim\cb$ and let
$$a_\vp=\sum_jv^{-2s_j}v^{d(\vp)}\in\ca$$ 
where $\r_{\cb!}\bbq=\op_j\bbq[-2s_j]$.

{\it Erratum to \cite{\GRA}.}
On page 202, line 1 of 16.8, replace "algebraic group $M$" by "algebraic 
group $M$ with a given Lie algebra homomorphism $\ph$ from $\fs$ to the Lie 
algebra of $M$".

On page 202, line 4 of 16.8, replace "Borel subgroups of $M$" by "Borel 
subgroups of the connected centralizer of $\ph(\fs)$ in $M$".

\subhead 11.2\endsubhead
The subset
$$\{a_\vp\i\tT_\boc\in\VV';\boc\in\un{\ovsc\EE},\vp\in\EE''\cap\boc\}\tag a$$
of $\VV'$ (see 10.14) is finite. Indeed, when $\boc$ is fixed, the subgroup 
$\fg^\ph_{\un0,\vp}$ of $\fg^\ph_{\un0}$ takes only finitely many values for 
$\vp$ in $\EE''\cap\boc$ hence $a_\vp$ takes only finitely many values.

Let $\VV'_\ca$ be the $\ca$-submodule of $\VV'$ generated by (a). Let 
$\VV_\ca$ be the image of $\VV'_\ca$ under the obvious linear map
$\VV'@>>>\VV$. The following result will be proved in 11.8.

(b) {\it $\VV_\ca$ is a free $\ca$-module such that the obvious $\QQ(v)$-linear
map \lb $\QQ(v)\ot_\ca\VV_\ca@>>>\VV$ is an isomorphism.}

\subhead 11.3\endsubhead
Let $\fh={}^\e\tfl^{\un\vp}$, $\fh^\ph=\fh\cap\fg^\ph$. We show:
$$\fg^\ph_{\un0,\vp}=\fh^\ph.\tag a$$
From 10.3(b) we have
$$\fh=\op_{N\in\ZZ,(\a,n)\in\car_{\un N};n=2N/\et,\la\vp:\a\ra=0}
(\fg^{\a,n}_{\un N}).$$ 
Using this and the definitions we have
$$\align&\{y\in\fh;[y,h]=0\}=\op_{N\in\ZZ,(\a,n)\in\car_{\un N};n=2N/\et=0,
\la\vp:\a\ra=0}(\fg^{\a,n}_{\un N})\\&
=\op_{(\a,n)\in\car_{\un0};n=0,\la\vp:\a\ra=0}(\fg^{\a,n}_{\un0})\endalign$$ 
and (a) follows.

\subhead 11.4\endsubhead
Recall that $\vp\in\EE''$. Let $\fp_*={}^{\e}\fp_*^{\un\vp}$, $\fu_*={}^{\e}\fu_*^{\un\vp}$.
We can find $\l'\in Y_Z$ such 
that $\la\l':\a\ra\ne0$ for any $i$ and any $(\a,n)\in\car^*_i$. 
Let $\vp'=\vp+\fra{1}{b}\l'\in\EE$. Let $\fp'_{*}={}^\e\fp^{\un\vp'}_*, 
\fu'_{*}={}^\e\fu^{\un\vp'}_*$ and $\fm'_{*}={}^\e\tfl^{\un\vp'}_*$. Assume that $b$ 
is sufficiently large; then $\vp'\in\EE'$ hence by 10.2 we have $\fm'_{*}=\fm_{*}$. The 
same argument as in 10.6(a) shows that $\fp'_N\sub\fp_N$ for all $N$. Since $b$ is large and 
$\vp'=\vp+\fra{1}{b}\l'$, we see that $\vp'$, $\vp$ are very close, so that 
$$\vp'\si\vp.\tag a$$ 
For $N\in\ZZ$ let $\fq_N$ be the image of $\fp'_N$ under the obvious 
projection $\fp_N@>>>\fh_N$. From 10.5(a),(b), we see that $\fq=\op_N\fq_N$ is
a parabolic subalgebra of $\fh$ and $\fm$ is a Levi subalgebra of $\fq$.
Moreover, if $u_N$ is the image of $\fu'_N$ 
under $\fp_N$ then $u=\op_Nu_N$ is the nilradical of $\fq$. 

From 10.3(b) we see that the $\ZZ$-grading of $\fh$ is $\et$-rigid and that
$e\in\ovsc\fh_{\et}$, 
so that $\ovsc\fm_\et\sub\ovsc\fh_{\et}$. 
Let $A_\vp\in\cq(\fh_\et)$ be the simple perverse sheaf on $\fh_\et$ such that
the support of $A_\vp$ is $\fh_\et$ and $A_\vp|_{\ovsc\fm_\et}$ is equal up to 
shift to $\tC|_{\ovsc\fm_\et}$.

Applying 1.8(b) and the transitivity formula 4.2(a) we deduce
$$\align& I_{\vp'}={}^\e\tInd_{\fp'_\et}^{\fg_\da}(\tC)=
{}^\e\tInd_{\fp_\et}^{\fg_\da}(\ind_{\fq_\et}^{\fh_\et}(\tC))
[\dim\fu'_0+\dim\fu'_\et-\dim\fu_0-\dim\fu_\et]\\&\cong
\op_j{}^\e\tInd_{\fp_\et}^{\fg_\da}(A_\vp)[-2s_j]
[\dim\fm_\et-\dim\fh_\et+\dim\fu'_0+\dim\fu'_\et-\dim\fu_0-\dim\fu_\et]\\&=
\op_j{}^\e\tInd_{\fp_\et}^{\fg_\da}(A_\vp)[-2s_j][d(\vp)],\tag b\endalign$$
where we have used the equality:
$$\dim\fm_\et-\dim\fh_\et+\dim\fu'_0+\dim\fu'_\et-\dim\fu_0-\dim\fu_\et=d(\vp).
\tag c$$
We now prove (c). Let $u'$ be the nilradical of the parabolic subalgebra
of $\fh$ that contains $\fm$ and is opposed to $\fq$. We have $u'=\op_Nu'_N$ 
where $u'_N=u'\cap\fh_N$. 
Now $u_0,u'_0$ are nilradicals of two opposite  
parabolic subalgebras of  $\fh_0$ hence $\dim u_0=\dim u'_0$.
We have $\dim\fh_\et=\dim u_\et+\dim u'_\et+\dim\fm_\et$,
$\dim\fu'_0-\dim\fu_0=\dim u_0$, $\dim\fu'_\et-\dim\fu_\et=\dim u_\et$. Hence the
left hand side of (b) equals $\dim u_0-\dim u'_\et=\dim u'_0-\dim u'_\et$.
Now $u'$ is normalized by $\fm$ hence by the Lie subalgebra $\fs$ of $\fm$
spanned by $e,h,f$. Note that $u'_0$ (resp. $u'_\et$) is the $0-$ (resp. $2-$)
eigenspace of $\ad(h):u'@>>>u'$. By the representation theory of $\fs$,
the map $\ad(e):u'_0@>>>u'_\et$ is surjective and its kernel is exactly the
space of $\fs$-invariants in $u'$ that is $u'\cap\fh^\ph$. We see that
$\dim u'_0-\dim u'_\et=\dim(u'\cap\fh^\ph)$.
Now $u'\cap\fh^\ph$ is the nilradical of a parabolic subalgebra of
$\fh^\ph$ with Levi subgroup $\fm\cap\fh^\ph=\fz$ (see 11.1) hence is the
nilradical of a Borel subalgebra of $\fh^\ph$ (which equals 
$\fg^\ph_{\un0,\vp}$ by 11.3(a)). By definition, the dimension of this 
nilradical is equal to $d(\vp)$. This proves (c) and hence also (b).

Now (b) implies the following equality in ${}^\x\ck(\fg_\da)$:
$$I_{\vp'}=a_\vp({}^\e\tInd_{\fp_\et}^{\fg_\da}(A_\vp)).$$
Using this and (a) we see that for any $\boc\in\un{\ovsc\EE}$ and any
$\vp\in\EE''\cap\boc$ we have
$$a_\vp\i\ti\g(\tT_\boc)={}^\e\tInd_{\fp_\et}^{\fg_\da}(A_\vp)
\in\QQ(v)\ot_\ca{}^\x\ck(\fg_\da).\tag d$$
Since ${}^\e\tInd_{\fp_\et}^{\fg_\da}(A_\vp)\in{}^\x\ck(\fg_\da)$, we see 
that 

(e) {\it for any $\boc\in\un{\ovsc\EE}$ and any $\vp\in\EE''\cap\boc$ we have 
$a_\vp\i\ti\g(\tT_\boc)\in{}^\x\ck(\fg_\da)$.}
\nl
The following result will be proved in 11.7.

(f) {\it The $\ca$-module ${}^\x\ck(\fg_\da)$ is generated by the elements
$a_\vp\i\ti\g(\tT_\boc)$ for various $\boc\in\un{\ovsc\EE}$ and 
$\vp\in\EE''\cap\boc$.}
\nl
From (f) we deduce:

(g) {\it The isomorphism $\g:\VV@>\si>>\QQ(v)\ot_\ca{}^\x\ck(\fg_\da)$ in 
10.19(a) restricts to an isomorphism of $\ca$-modules
$\g_\ca:\VV_\ca@>\si>>{}^\x\ck(\fg_\da)$.}

\subhead 11.5\endsubhead
Let $\fp_*$ be an $\e$-spiral with a splitting $\fh_*$ and with nilradical
$\fu_*$; let $A\in\cq(\fh_\et)$ be a simple perverse sheaf. We show:

(a) {\it Let $\cx$ be the collection of all $B\in\fB$ such that some shift of 
$B$ is a direct summand of $\tInd_{\fp_\et}^{\fg_\da}(A)$. Then the map 
$\cx@>>>\un\fT_\et$, $B\m\ps(B)$ is constant.}
\nl
We can find a parabolic subalgebra $\fq$ of $\fh$ and a Levi subalgebra $\fm'$
of $\fq$ such that $\fq=\op_N\fq_N,\fm'=\op_N\fm'_N$ where 
$\fq_N=\fq\cap\fh_N$, $\fm'_N=\fm'\cap\fh_N$ and a cuspidal perverse sheaf 
$C$ in $\cq(\fm_\et)$ such 
that some shift of $A$ is a direct summand of $\ind_{\fq_\et}^{\fh_\et}(C)$. Let
$M'=e^{\fm'},M'_0=e^{\fm'_0}$. Setting $\fp'_N=\fu_N+\fq_N$ 
for any $N\in\ZZ$,  
we see from 2.8(a) that $\fp'_*$ is an $\e$-spiral and from 2.8(b) that 
$\fm'_*$ is a splitting of $\fp'_*$. We see that 
$(M',M'_0,\fm',\fm'_*,C)\in\fT_\et$. Let 
$\x'$ be the element of $\un\fT_\et$ determined by $(M',M'_0,\fm',\fm'_*,C)$. 
If $B\in\cx$ then (by 4.2)
some shift of $B$ is a direct summand of $\tInd_{\fp'_\et}^{\fg_\da}(C)$ hence 
$\ps(B)=\x'$. This proves (a).

\mpb

We say that $\tInd_{\fp_\et}^{\fg_\da}(A)$ (as in (a)) has type 
$\x'\in\un\fT_\et$ if $\ps(B)=\x'$ for any $B\in\cx$.

\subhead 11.6\endsubhead
Recall from 8.4(a) that the $\ca$-module $\ck(\fg_\da)$ is generated by the
classes of $\e$-quasi-monomial objects of $\cq(\fg_\da)$. Using this and 
11.5(a), we deduce that in the direct sum decomposition
$\ck(\fg_\da)=\op_{\x\in\un\fT_\et}{}^\x\ck(\fg_\da)$ (see 6.7), any summand
${}^\x\ck(\fg_\da)$ is generated as an $\ca$-module by the classes of 
$\et$-quasi-monomial objects in $\cq(\fg_\da)$ of type $\x$. 

\subhead 11.7\endsubhead
We prove 11.4(f). (Thus we are again in the setup of 11.1.) 
Using 11.6, we see that it is enough to show that if $A'$ is
an $\et$-quasi-monomial object in $\cq(\fg_\da)$ of type $\x$, then the class of
$A'$ in $\ck(\fg_\da)$ is of the form $a_\vp\i\ti\g(\tT_\boc)$ for some 
$\boc\in\un{\ovsc\EE}$ and $\vp\in\EE''\cap\boc$. We can find

(a) $\fp_*,\fh_*,A,\fq_*,\fp'_*$, $(M',M'_0,\fm',\fm'_*,C)\in\fT_\et$ 
(representing $\x$)
\nl
as in 11.5 such that $A'=\tInd_{\fp_\et}^{\fg_\da}(A)$; moreover, we can assume 
that the $\ZZ$-grading of $\fh=\op_N\fh_N$ is $\et$-rigid and 
$\ovsc{\fm'_\et}\sub\ovsc{\fh_\et}$. Replacing the data (a) by a 
$G_{\un0}$-conjugate we can assume in addition that $(M',M'_0,\fm',\fm'_*,C)$ 
is equal to $(M,M_0,\fm,\fm_*,\tC)$ in 10.1.

Let 
$H=e^{\fh}$. Since $\fh_*$ is $\et$-rigid, there exists $\io'\in Y_H$ such that

(i) ${}^{\io'}_k=\fh_{k\et/2}$ if $k\in\ZZ$, $k\et/2$,  ${}^{\io'}_k=0$ if $k\in\ZZ,k\et/2\n\ZZ$
and

(ii) $\io'=\io_{\ph'}$ for some $\ph'=(e',h',f')\in J^H$ such that 
$e'\in\ovsc\fh_\et$, $h'\in\fh_0$, $f'\in\fh_{-\et}$.

Since $\ovsc\fm_\et\sub\ovsc\fh_\et$ and $e\in\ovsc\fm_\et$, we see that $e,e'$ 
are in the same $M_0$-orbit. Hence we can find $g\in M_0$ such that $\Ad(g)$ 
conjugates $e',f',h',\io'$ to $e,f,h,\io$. Applying $\Ad(g)$ (which preserves 
$\fh_k$) to (i) and (ii) we see that we can assume that $\io'=\io$, $\ph'=\ph$.

Recall from 10.3 that $\tfl^\ph_N={}_{2N/\et}^\io\fg_{\un N}$ for $N\in\ZZ$ such that $2N/\et\in\ZZ$. 
Using (i) with $\io'=\io$ we see that $\fh_N\sub\tfl^\ph_N$ for any $N\in\ZZ$ such that $2N/\et\in\ZZ$. 
Using 10.4(c),(d) we see that for some $\vp\in\EE''$ we have
$\fp_*={}^\e\fp^{\un\vp}_*$, $\fh_*={}^\e\tfl^{\un\vp}_*$. Using now 11.4(d), 
we see that $A'=a_\vp\i\ti\g(\tT_\boc)$, where $\boc\in\un{\ovsc\EE}$ 
contains $\vp$. This completes the proof of 11.4(f), hence that of 11.4(g).

\subhead 11.8\endsubhead
We can now prove 11.2(b). Using 11.4(g), 11.2(b) is reduced to the following
obvious statement: ${}^\x\ck(\fg_\da)$ is a free $\ca$-module.

\subhead 11.9\endsubhead
We define a $\QQ$-linear map 
$$\bar{}:\QQ(v)\ot_\ca\ck(\fg_\da)@>>>\QQ(v)\ot_\ca\ck(\fg_\da)$$ by 
$\ov{fB}=\baf B$ for any $f\in\QQ(v)$ and any $B\in\fB$ (see 10.13); here 
$\baf$ is as in 0.12. This restricts to a $\QQ$-linear map 
$$\bar{}:\QQ(v)\ot_\ca{}^\x\ck(\fg_\da)@>>>\QQ(v)\ot_\ca{}^\x\ck(\fg_\da)$$ 
and to a $\ZZ$-linear map ${}^\x\ck(\fg_\da)@>>>{}^\x\ck(\fg_\da)$. We show:

(a) $\ov{I_\vp}=I_\vp$ for any $\vp\in\EE'$.
\nl
In ${}^\x\cq(\fg_\da)$ we have $I_\vp=\sum_{B\in{}^\x\fB}f_BB$
where $f_B\in\ca$. It is enough to prove that $\baf_B=f_B$ for all $B$.

We set $\fp_*={}^\e\fp^{\un\vp}_*$.
Let $\s$ be an automorphism of order $2$ of $\bbq$ such that $\s(z)=z\i$ for
any root of $1$ in $\bbq$. Applying $\s$ to $K\in\cq(\fm_\et)$ (resp.
$K\in\cq(\fg_\da)$) we obtain $K^\s\in\cq(\fm_\et)$ (resp. $K^\s\in\cq(\fg_\da)$).
Note that $K\m K^\s$ commutes with shifts; moreover we have
$$({}^\e\tInd_{\fp_\et}^{\fg_\da}(\tC))^\s={}^\e\tInd_{\fp_\et}^{\fg_\da}(\tC^\s).
\tag b$$
Moreover, if $K$ is a simple perverse sheaf in $\cq(\fm_\et)$ or in 
$\cq(\fg_\da)$, we have
$$K^\s\cong D(K),\tag c$$
since $K$ restricted to an open dense subset of it support is a local system 
with finite monodromy.
By (b),(c) and 4.1(d) we have
$$D({}^\e\tInd_{\fp_\et}^{\fg_\da}(\tC))={}^\e\tInd_{\fp_\et}^{\fg_\da}(D(\tC))   
={}^\e\tInd_{\fp_\et}^{\fg_\da}(\tC^\s)=({}^\e\tInd_{\fp_\et}^{\fg_\da}(\tC))^\s,$$
hence 
$$(D(I_\vp))^\s=I_\vp$$
and $\sum_B\baf_BD(B^\s)=\sum_B f_BB$. Using this and (c) for $K=B$, we see
that $\baf_B=f_B$ for all $B$; this proves (a).

\subhead 11.10\endsubhead
We show:

(a) {\it $\bar{}:\VV'@>>>\VV'$ (see 10.21) restricts to an involution 
$\VV'_\ca@>>>\VV'_\ca$ and $\bar{}:\VV@>>>\VV$ (see 10.21) restricts to an 
involution $\VV_\ca@>>>\VV_\ca$ (these restrictions are denoted again by 
$\bar{}$.}
\nl
It is enough to note that for any $\vp\in\EE''$ we have $\ov{a_\vp}=a_\vp$.

\mpb

We show:

(b) {\it $\ti\g:\VV'@>>>\QQ(v)\ot_\ca{}^\x\ck(\fg_\da)$ is compatible with the 
maps $\bar{}$ in the two sides.}
\nl
This follows from 11.9(a).

\mpb

We have the following result.
\proclaim{Proposition 11.11} (a) Let $\BB'$ be the set of all $b\in\VV_\ca$ 
such that $\bar b=b$ and $(b:b)\in1+v\ZZ[[v]]$. Then $\BB'$ is a signed basis 
of the $\ca$-module $\VV_\ca$ (that is, the union of a basis with $(-1)$ times
that basis).

(b) For $b\in\BB'$ we have $(b:b)\in1+v\NN[[v]]$. 

(c) There is a unique $\ca$-basis $\BB$ of $\VV_\ca$ such that $\BB\sub\BB'$
and for any $\boc\in\un{\ovsc\EE}$ and any $\vp\in\EE''\cap\boc$, the image of
$a_\vp\i\tT_\boc$ in $\VV_\ca$ is an $\NN[v,v\i]$-linear combination of 
elements in $\BB$.
\endproclaim
It is enough to prove the analogous statements where $\VV_\ca$ is identified
via $\g$ with ${}^\x\ck(\fg_\da)$ with $(:)$ as in 4.4(c) and with $\bar{}$ as 
in 11.9 (we use 10.19(b), 11.10(b)). 
Let ${}^\x\fB=\{B_1,B_2,\do,B_r\}$ (see 10.13). From 0.12 we have 
$(B_j:B_{j'})\in\d_{j,j'}+h_{j,j'}$ where $h_{j,j'}\in v\NN[[v]]$ for all 
$j,j'$ in $[1,r]$. From the definition 
(see 11.9) we have $\ov{B_j}=B_j$ for $j=1,\do,r$. Now let
$b\in{}^\x\ck(\fg_\da)$ be such that $\bar b=b$ and $(b:b)\in1+v\ZZ[[v]]$. To
prove (a), it is enough to show that $b=\pm B_j$ for some $j$. We can write 
$b=\sum_{j=1}^rf_jB_j$ where $f_j\in\ca$ satisfy $\baf_j=f_j$ and
$\sum_{j,j'\in[1,r]}\baf_jf_{j'}(\d_{j,j'}+h_{j,j'})\in1+v\ZZ[[v]]$ 
hence $\sum_{j,j'\in[1,r]}f_jf_{j'}(\d_{j,j'}+h_{j,j'})\in1+v\ZZ[[v]]$. 
We can find
$c\in\ZZ$ such that $f_j=f_{j,c}v^c\mod v^{c+1}\ZZ[v]$ where 
$f_{j,c}\in\ZZ$ for all $j$ and $f_{j,c}\ne0$ for some $j$. We have
$\sum_{j\in[1,r]}f_{j,c}^2v^{2c}+v^{2c+1}f'=1+vf''$ where $f',f''\in\ZZ[v]$. 
Moreover, $\sum_{j\in[1,r]}f_{j,c}^2>0$. It follows that $c=0$ and 
$\sum_{j\in[1,r]}f_{j,0}^2=1$ so that there exists $j_0\in[1,r]$ such that 
$f_{j_0,0}=\pm1$ and $f_{j,0}=0$ for $j\ne j_0$. We have 
$f_j=\pm\d_{j,j_0}\mod v\ZZ[v]$ for all $j$. Since $\baf_j=f_j$ we deduce 
that $f_j=\pm\d_{j,j_0}$ for all $j$. Thus $b=\pm B_{j_0}$. This completes the
proof of (a). At the same time we have proved (b). Clearly, 
$\{B_1,B_2,\do,B_r\}$ has the positivity property in (c) (with $\VV_\ca$ 
identified with ${}^\x\ck(\fg_\da)$ and with $a_\vp\i\tT_\boc$ identified with 
$\g(a_\vp\i\tT_\boc)$). Since any $B_j$ appears with $>0$ coefficient in some
$\g(a_\vp\i\tT_\boc)$, we see that $\{B_1,B_2,\do,B_r\}$ is the only basis
contained in $\{\pm B_1,\pm B_2,\do,\pm B_r\}$ with the positivity property in
(c). This completes the proof of the proposition.

\subhead 11.12\endsubhead
From the proof of 11.11 we see that 
$\g:\VV@>\si>>\QQ(v)\ot_\ca{}^\x\ck(\fg_\da)$ (see 10.19(a)) restricts to a 
bijection
$$\BB@>\si>>{}^\x\fB.\tag a$$
For any $G_{\un0}$-orbit $\co$ in $\fg_\da^{nil}$ let $\BB_\co$ be the set of 
all $b\in\BB$ such that $\g(b)\in\fB_\co$ (see 10.13). 
We have a partition $\BB=\sqc_\co\BB_\co$ where $\co$ runs over the
$G_{\un0}$-orbits in $\fg_\da^{nil}$.

\subhead 11.13\endsubhead
We consider the setup of 10.22. We show:

(a) {\it the $\ca$-submodule $\VV_\ca$ of $\VV$ defined in 11.2 in terms of $\dot\x$ is the
same as that defined in terms of $\dot\x_1$.}
\nl
It is not clear how to prove this using the definition in 11.2 since $\EE''$ defined in terms
of $\dot\x$ (see 10.3) is not necessarily the same as that defined in terms of $\dot\x_1$.
Instead we will argue indirectly.
Using 11.4(g) it is enough to show 

(b) {\it the isomorphism 
$\g:\VV@>\si>>\QQ(v)\ot_\ca{}^\x\ck(\fg_\da)$ defined in 10.19 in terms of
$\dot\x$ is equal to the analogous isomorphism defined in terms of $\dot\x_1$.}
\nl
Thus it is enough to show that if $\vp\in\EE'\cap\EE'_1$ then
$\tI_\vp$ defined in 10.2 in terms of $\dot\x$ is the same as that defined in terms of
$\dot\x_1$.
Using the definitions we see that it is enough to show that
$${}^{\doet}\fp^{(|\et|/2)(\vp+\io)}_\et={}^{\doet_1}\fp^{(|\et_1|/2)(\vp+\io)}_{\et_1},$$
$${}^{\doet}\fp^{(|\et|/2)(\vp+\io)}_0={}^{\doet_1}\fp^{(|\et_1|/2)(\vp+\io)}_0.$$
or that
$$\op_{\k\in\QQ;\k\ge|\et|}({}^{(|\et|/2)(\vp+\io)}_\k\fg_\d)=
\op_{\k\in\QQ;\k\ge|\et_1|}({}^{(|\et_1|/2)(\vp+\io)}_\k\fg_\d),$$
$$\op_{\k\in\QQ;\k\ge0}({}^{(|\et|/2)(\vp+\io)}_\k\fg_\d)=
\op_{\k\in\QQ;\k\ge0}({}^{(|\et_1|/2)(\vp+\io)}_\k\fg_\d),$$
or that
$$\op_{\k\in\QQ;\k/|\et|\ge1}({}^{(1/2)(\vp+\io)}_{\k/|\et|}\fg_\d)=
\op_{\k\in\QQ;\k/|\et_1|\ge1}({}^{(1/2)(\vp+\io)}_{\k/|\et_1|}\fg_\d),$$
$$\op_{\k\in\QQ;\k/|\et|\ge0}({}^{(1/2)(\vp+\io)}_{\k/|\et|}\fg_\d)=
\op_{\k\in\QQ;\k/|\et_1|\ge0}({}^{(1/2)(\vp+\io)}_{\k/|\et_1|}\fg_\d),$$
or, setting $\k'=\k/|\et|$, $\k''=\k/|\et_1|$, that
$$\op_{\k'\in\QQ;\k'\ge1}({}^{(1/2)(\vp+\io)}_{\k'}\fg_\d)=
\op_{\k''\in\QQ;\k''\ge1}({}^{(1/2)(\vp+\io)}_{\k''}\fg_\d),$$
$$\op_{\k'\in\QQ;\k'\ge0}({}^{(1/2)(\vp+\io)}_{\k'}\fg_\d)=
\op_{\k''\in\QQ;\k''\ge0}({}^{(1/2)(\vp+\io)}_{\k''}\fg_\d),$$
which are obvious. This proves (a).

\mpb

Using (b) and 11.12(b) we see that the basis $\BB$ of $\VV_\ca$
defined in 11.11 in terms of 
$\dot\x$ is the same as that defined in terms of $\dot\x_1$.

\head 12. Purity properties\endhead
In this section we show that for any irreducible local system $\cl$
on a $G_{\un0}$-orbit in $\fg^{nil}_\da$ the cohomology sheaves of 
$\cl^\sha\in\cd(\fg_\da)$ satisfy a strong purity property. This generalizes
the analogous result in the $\ZZ$-graded case in \cite{\GRA}.

\subhead 12.1\endsubhead
In this section we assume that $p>0$ and that $\kk$ is an algebraic closure of
a finite field $\FF_q$ with $q$ elements (here $q$ is a power of $p$). 
Replacing $q$ by larger powers of $p$ if necessary,  
we can assume that $m$ divides $q-1$ and that we 
can find an $\FF_q$-rational structure on $G$ with Frobenius map $F:G@>>>G$ 
such that $\vt:G@>>>G$ (see 0.5) commutes with
$F:G@>>>G$. Then $G_{\un0}$ is defined over $\FF_q$ and
$\fg$ inherits from $G$ an $\FF_q$-rational structure
with Frobenius map $F:\fg@>>>\fg$ satisfying $F(\fg_i)=\fg_i$ for all $i$.
Again by replacing $q$ by larger powers of $p$ if necessary, 
we may assume that all $G_{\un0}$-orbits in $\fg_\da^{nil}$ are defined over 
$\FF_q$ and that for any irreducible $G_{\un0}$-equivariant
local system $\cl$ on a $G_{\un0}$-orbit $\co$ in $\fg_\da^{nil}$ we have 
$F^*\cl\cong\cl$. 

We now fix a $G_{\un0}$-orbit $\co$ in $\fg_\da^{nil}$ with closure $\bco$ and 
an irreducible $G_{\un0}$-equivariant local system 
$\cl$ on $\co$. We fix an isomorphism $\tF:F^*\cl@>>>\cl$ which induces the
identity map on the stalk of $\cl$ at some point of $\co^F$. Then $\tF$ induces
an isomorphism (denoted again by $\tF$) $F^*\cl^\sha@>>>\cl^\sha$. 
Given a finite dimensional $\bbq$-vector space $V$ with an endomorphism 
$\tF:V@>>>V$, we say that $\tF:V@>>>V$ is $a$-pure (for an integer $a$) if 
the eigenvalues of $\tF$ are algebraic numbers all of whose complex conjugates
have absolute value $q^{a/2}$. Sometimes we will just say that $V$ is $a$-pure
(where this is understood to refer to $\tF$).

\subhead 12.2\endsubhead
We show:

(a) {\it For any $x\in\bco^F$ and any $a\in\ZZ$, the induced linear map 
$\tF:\ch^a_x(\cl^\sha)@>>>\ch^a_x(\cl^\sha)$ is $a$-pure.}
\nl
Using 2.3(b), we can find $\ph=(e,h,f)\in J_\da(x)$ such that $h,f$ 
are $\FF_q$-rational. Let $\io=\io_\ph\in Y_G$ be as in 1.1. Let $\fz(f)$ be the
centralizer of $f$ in $\fg$ and let $\Si=e+\fz(f)\sub\fg$. According to
Slodowy, 

(b) {\it the affine space $\Si$ is a transversal slice at $x$ to the $G$-orbit
of $e$ in $\fg$ and the $\kk^*$-action $t\m t^{-2}\Ad(\io(t))$ on $\fg$ keeps 
$e$ fixed, leaves $\Si$ stable and defines a contraction of $\Si$ to $x$.}
\nl
Let $\ti\Si=\Si\cap\fg_\da=e+(\fz(f)\cap\fg_\da)$. Then

(c) {\it $\ti\Si$ is a transversal slice to the $G_{\un0}$-orbit of $e$ in 
$\fg_\da$; the $\kk^*$-action in (b) leaves stable $\ti\Si$ and is a 
contraction of $\ti\Si$ to $e$.}
\nl
(We have a direct sum decomposition $\fg_\da=(\fz(f)\cap\fg_\da)\op[e,\fg_0]$
obtained by taking the $\z^\da$-eigenspace of $\th:\fg@>>>\fg$ in the two sides
of the direct sum decomposition $\fg=\fz(f)\op[e,\fg]$; note that both
$\fz(f)$ and $[e,\fg]$ are $\th$-stable.)

Let $\cl'$ be the restriction of $\cl$ to $\co\cap\ti\Si$ (a smooth 
irreducible subvariety of $\ti\Si$). Note that $\co\cap\ti\Si$,
$\bco\cap\ti\Si$ are stable under the $\kk^*$-action in (c) and $\cl'$ is
equivariant for that action. By (c), we have 
$\ch_x(\cl^\sha)=\ch_x(\cl'{}^\sha)$. It remains to show that
$\ch_x(\cl'{}^\sha)$ is $a$-pure. This can be reduced to Deligne's hard
Lefschetz theorem by an argument in Lemma 4.5(b) in \cite{\KLL} applied to
$\bco\cap\ti\Si\sub\ti\Si$ with the contraction in (c) and to $\cl'$ instead of
$\bbq$. Note that in \cite{\KLL, 4.5(b)} an inductive purity assumption was
made which is in fact unnecessary, by Gabber purity theorem. This completes the
proof of (a).

{\it Erratum to \cite{\GRA}. On page 209, line -5, replace 
$t^{-n}\Ad(\io'(t))$ by $t^{-2}\Ad(\io'(t))$.}

\head 13. An inner product\endhead
This section is an adaptation of \cite{\GRA,\S19} from the $\ZZ$-graded case
to the $\ZZ/m$-graded case; we express  
the matrix whose entries are the 
values of the $(:)$-pairing at two elements of $\fB$ as a product of three 
matrices. As an application we show (see 13.8(a)) that if
$(\co,\cl),(\co',\cl')$ in $\ci(\fg_\da)$ are such that some
cohomology sheaf of $\cl^\sha|_{\co'}$ contains $\cl'$ then
$(\co,\cl),(\co',\cl')$ are in the same block. Another application (to odd
vanishing) is given in Section 14.

\subhead 13.1\endsubhead
We fix $(\co,\cl),(\co',\cl')$ in $\ci(\fg_\da)$ and we form $A=\cl^\sha$,
$A'=\cl'{}^\sha$ in $\cd(\fg_\da)$. We want to compute $d_j(\fg_\da;A,A')$ (see 
0.12) for a fixed $j\in\ZZ$.
We can arrange the $G_{\un0}$-orbits in $\fg_\da^{nil}$ in an order
$\co_0,\co_1,\do,\co_\b$ such that $\co_{\le s}=\co_0\cup\co_1\cup\do\cup\co_s$
is closed in $\fg_\da$ for $s=0,1,\do,\b$. We choose a smooth irreducible
variety $\G$ with a free action of $G_{\un0}$ such that $H^r(\G,\bbq)=0$ for
$r=1,2,\do,\mm$ where $\mm$ is a large integer (compared to $j$). We assume
that $\DD:=\dim\G$ is large (compared to $j$). We have 
$H^{2\DD-i}_c(\G,\bbq)=0$ for $i=1,2,\do,\mm$.
We form $X=G_{\un0}\bsl(\G\T\fg_\da)$. 
Let $\tcl,\tcl'$ be the local systems on the smooth subvarieties 
$G_{\un0}\bsl(\G\T\co)$, $G_{\un0}\bsl(\G\T\co')$ of $X$ defined by $\cl,\cl'$
and let $\tA=\tcl^\sha$, $\tA'=\tcl'{}^\sha$ be the corresponding intersection
cohomology complexes in $\cd(X)$. Then $\tA\ot\tA'\in\cd(X)$ is well defined;
its restriction to various subvarieties of $X$ will be denoted by the same
symbol. For $s=0,1,\do,\b$ we form
$X_s=G_{\un0}\bsl(\G\T\co_s)$,  $X_{\le s}=G_{\un0}\bsl(\G\T\co_{\le s})$.
We set $X_{\le-1}=\emp$. The partition of $X_{\le s}$ into $X_{\le s-1}$ and
$X_s$ (for $s=0,1,\do,\b$) gives rise to a long exact sequence
$$\align&\do@>\x_a>>H^a_c(X_s,\tA\ot\tA')@>>>H^a_c(X_{\le s},\tA\ot\tA')@>>>\\&
H^a_c(X_{\le s-1},\tA\ot\tA')@>\x_{a+1}>>H^{a+1}_c(X_s,\tA\ot\tA')@>>>\do
\tag a\endalign$$

\subhead 13.2\endsubhead
In the following proposition we encounter two kinds of integers; some like $j$,
$m$, $\dim G_{\un0}$, $\dim\fg_\da$, $\b$ are regarded as "small" (they belong 
to a fixed finite set of integers), others like $2\DD$ and $\mm$ are regarded 
as very large (we are free to choose them so). We will also encounter integers
$a$ such thay $2\DD-a$ is a "small" integer (we then write $a\si2\DD$).

\proclaim{Proposition 13.3} (a) Assume that $\kk$ is as in 12.1. Then
$H^a_c(X_s,\tA\ot\tA')$ is $a$-pure (see 12.1) for $s=0,1,\do,\b$ and 
$a\si2\DD$.

(b) We choose $x_s\in\co_s$. If $a\si2\DD$ then 
$$\align&\dim H^a_c(X_s,\tA\ot\tA')\\&=\sum_{r+r_1+r_2=a}
\dim(H^r_c(G_{\un0}(x_s)^0\bsl\G,\bbq)
\ot\ch^{r_1}_{x_s}A\ot\ch^{r_2}_{x_s}A')^{G_{\un0}(x_s)/G_{\un0}(x_s)^0}
\endalign$$
where the upper script refers to invariants under the finite group
$G_{\un0}(x_s)/G_{\un0}(x_s)^0$.

(c) Assume that $\kk$ is as in 12.1. Then $H^a_c(X_{\le s},\tA\ot\tA')$ is 
$a$-pure (see 12.1) for $s=0,1,\do,\b$ and $a\si2\DD$.

(d) The exact sequence 13.1(a) gives rise to short exact sequences
$$0@>>>H^a_c(X_s,\tA\ot\tA')@>>>H^a_c(X_{\le s},\tA\ot\tA')@>>>
H^a_c(X_{\le s-1},\tA\ot\tA')@>>>0$$
for $s=0,1,\do,\b$ and $a\si2\DD$.

(e) For $a\si2\DD$ we have 
$\dim H^a_c(X,\tA\ot\tA')=\sum_{s=0}^\b\dim H^a_c(X_s,\tA\ot\tA')$.

(f) For $a\si2\DD$ we have 
$$\align&\dim H^a_c(X,\tA\ot\tA')\\&=\sum_{s=0}^\b\sum_{r+r_1+r_2=a}
\dim(H^r_c(G_{\un0}(x_s)^0\bsl\G,\bbq)\ot\ch^{r_1}_{x_s}A\ot
\ch^{r_2}_{x_s}A')^{G_{\un0}(x_s)/G_{\un0}(x_s)^0}.\endalign$$
\endproclaim
The proof is almost a copy of the proof of \cite{\GRA, 19.4}.
By general principles we can assume that $\kk$ is as in 12.1. We shall use
Deligne's theory of weights. We first prove (a) and (b). We write $x$ instead
of $x_s$. We may assume that $x$ is an $\FF_q$-rational point. We have a
natural spectral sequence
$$E^{r,r'}_2=H^r_c(X_s,\ch^{r'}(\tA\ot\tA'))\imp H^{r+r'}_c(X_s,\tA\ot\tA').
\tag g$$
We show that 
$$E^{r,r'}_2 \text{ is $(r+r')$-pure if }r+r'\si2\DD.\tag h$$
We have $X_s=G_{\un0}(x)\bsl\G$ and
$$E^{r,r'}_2=(H^r_c(G_{\un0}(x)^0\bsl\G,\bbq)\ot
\ch^{r'}_x(A\ot A'))^{G_{\un0}(x)/G_{\un0}(x)^0}.$$
Here $A\ot A'\in\cd(\fg_\da)$. We may assume that $r'$ is "small" (otherwise,
$E^{r,r'}_2=0$). We then have $r\si2\DD$. Now
$$H^{r'}_x(A\ot A')=\op_{r_1+r_2=r'}\ch^{r_1}_x(A)\ot\ch^{r_2}_x(A')$$
is $r'$-pure by 12.2(a). Moreover, it is well known that 
$H^r_c(G_{\un0}(x)^0\bsl\G,\bbq)$ is $r$-pure for $r\si2\DD$ and (h) follows.

From (h) it follows that $E^{r,r'}_{\iy}$ of the spectral sequence (g) is
$(r+r')$-pure if $r+r'\si2\DD$ and (a) follows. From (h) it also follows that
$E^{r,r'}_2=E^{r,r'}_{\iy}$ if $r+r'\si2\DD$ (many differentials must be zero
since they respect weights) and (b) follows.

Now (c) follows from (a) using the exact sequence 13.2(a) and induction on $s$.
From (c) and (a) we see that the homomorphism $\x_{a+1}$ in 13.2(a) is between
pure spaces of different weight. Since $\x_{a+1}$ preserve weights, it must be
$0$. Similarly $\x_a=0$ in 13.2(a) hence (d) holds. Now (e) follows from (d)
since the support of $\tA\ot\tA'$ is contained in $X_{\le\b}$; (f) follows from
(b),(e). This completes the proof of the proposition.

\subhead 13.5\endsubhead
Using the definitions we see that 13.4(f) implies:
$$d_j(\fg_\da;A,A')=\sum_s\sum_{-r_0+r_1+r_2=j-2p_s}
\dim(H_{r_0}^{G_{\un0}(x_s)^0}(.)
\ot\ch^{r_1}_{x_s}A\ot\ch^{r_2}_{x_s}A')^{G_{\un0}(x_s)/G_{\un0}(x_s)^0}\tag a
$$
where 
$$p_s=\dim G_{\un0}-\dim G_{\un0}(x_s)=\dim\co_s\tag b$$
and $H_{r_0}^{G_{\un0}(x_s)^0}(.)$ denotes equivariant homology of a point.
(See \cite{\CUS} for the definition of equivariant homology.)

\subhead 13.6\endsubhead
Given $(\co,\cl),(\ti\co,\tcl)$ in $\ci(\fg_\da)$ we define 
$\mu(\tcl,\cl)\in\ZZ[v\i]$ by
$$\mu(\tcl,\cl)=\sum_a\mu(a;\tcl,\cl)v^{-a}\tag a$$
where $\mu(a;\tcl,\cl)$ is the number of times $\tcl$ appears in a
decomposition of the local system $\ch^a(\cl^\sha)|_{\ti\co}$ as a direct sum of
irreducible local systems. Note that $\mu(\tcl,\cl)$ is zero unless $\ti\co$ is
contained in the closure of $\co$.

If $\ce$ is an irreducible $G_{\un0}$-equivariant local system on $\co_s$, we
denote by $\r_\ce$ the irreducible $G_{\un0}(x_s)/G_{\un0}(x_s)^0$-module
corresponding to $\ce$. With this notation we can rewrite 13.5(a) as follows:
$$\align&d_j(\fg_\da;A,A')=\sum_s\sum_{-r_0+r_1+r_2=j-2p_s}\\&\sum_{\ce,\ce'}
\mu(r_1;\ce,\cl)\mu(r_2;\ce',\cl')\dim(H_{r_0}^{G_{\un0}(x_s)^0}(.)\ot
\r_\ce\ot\r_{\ce'})^{G_{\un0}(x_s)/G_{\un0}(x_s)^0}\endalign$$
where $\ce,\ce'$ run over the isomorphism classes of irreducible 
$G_{\un0}$-equivariant local systems on $\co_s$. This may be written in terms 
of power series in $\QQ((v))$ as follows.
$$\{A,A'\}=\sum_{s=0}^\b\sum_{\ce,\ce'}\mu(\ce,\cl)\Xi(\ce,\ce')\mu(\ce',\cl')
\tag b$$
where
$$\Xi(\ce,\ce')=\sum_{r_0\ge0}\dim(H_{r_0}^{G_{\un0}(x_s)^0}(.)\ot
\r_\ce\ot\r_{\ce'})^{G_{\un0}(x_s)/G_{\un0}(x_s)^0}v^{r_0-2p_s}\in\QQ((v)).
\tag c$$

\subhead 13.7\endsubhead
Let $B_1,B_2\in\fB$. We write $B_1=\cl_1^\sha[\dim\co_1]\in\cd(\fg_\da)$, 
$B_2=\cl_2^\sha[\dim\co_2]\in\cd(\fg_\da)$ with $(\co_1,\cl_1),(\co_2,\cl_2)$ 
in $\ci(\fg_\da)$. We set

$P_{B_2,B_1}=\mu(\cl_2,\cl_1)\in\NN[v\i]$, (see 13.6);

$\tP_{B_2,B_1}=P_{D(B_2),D(B_1)}\in\NN[v\i]$;

$\L_{B_2,B_1}=\Xi(\cl_2,\cl_1)\in\QQ((v))$, (see 13.6) if $\co_1=\co_2$;

$\L_{B_2,B_1}=0$ if $\co_1\ne\co_2$;

$\ti\L_{B_2,B_1}=\L_{B_2,D(B_1)}\in\QQ((v))$.
\nl
Note that 

(a) $P_{B_2,B_1}\ne0\imp\dim\co_2\le\dim\co_1$;
$\tP_{B_2,B_1}\ne0\imp\dim\co_2\le\dim\co_1$;

(b) if $\dim\co_2=\dim\co_1$, then $P_{B_2,B_1}=\d_{B_2,B_1}$,
$\tP_{B_2,B_1}=\d_{B_2,B_1}$;

(c) $\ti\L_{B_2,B_1}=0$ if $\co_1\ne\co_2$.
\nl
Then 13.6(b) can be rewritten as follows:
$$\{A,A'\}=\sum_{B_1\in\fB,B_2\in\fB}P_{B_1,B}\L_{B_1,B_2}P_{B_2,B''}$$
where $B=A[\dim\co]$, $B''=A'[\dim\co']$, or as
$$\{A,A'\}=\sum_{B_1\in\fB,B_2\in\fB}P_{B_1,B}\L_{B_1,D(B_2)}P_{D(B_2),D(B')}$$
where $B=A[\dim\co]$, $D(B')=A'[\dim\co']$, or as
$$\{A,A'\}=\sum_{B_1\in\fB,B_2\in\fB}P_{B_1,B}\ti\L_{B_1,B_2}\tP_{B_2,B'}$$
where $B=A[\dim\co]$, $D(B')=A'[\dim\co']$. We have 
$$\{A,A'\}=v^{-\k(B)-\k(B')}\{B,D(B')\}=v^{-\k(B)-\k(B')}(B:B'),$$
hence
$$v^{-\k(B)-\k(B')}(B:B')=
\sum_{B_1\in\fB,B_2\in\fB}P_{B_1,B}\ti\L_{B_1,B_2}\tP_{B_2,B'}\tag d$$
for any $B,B'$ in $\fB$. (Here $\k$ is as in 10.13.) We show:

(e) {\it the following three matrices with entries in $\QQ((v))$ (indexed by
$\fB\T\fB$) are invertible:

(i) the matrix $((B:B'))$;

(ii) the matrix $\cm:=(v^{-\k(B)-\k(B')}(B:B'))$;

(iii) the matrix $\ct:=(\ti\L_{B,B'})$.}
\nl
The matrix in (i) is invertible since $(B:B')\in\d_{B,B'}+v\NN[[v]]$ for all 
$B,B'$, see 0.12. This implies immediately that $\cm$ is invertible. Now, by 
(d), we have $\cm=\cs\ct\cs'$ where $\cs$ (resp. $\cs'$) is the matrix indexed
by $\fB\T\fB$ whose $(B,B')$-entry is $P_{B',B}$ (resp. $\tP_{B,B'}$). Since
$\cm$ is invertible, it follows that $\ct$ is invertible. This proves (e).

\subhead 13.8\endsubhead
We show:

(a) {\it If $B,B'$ in $\fB$ satisfy $\psi(B)\ne\psi(B')$ ($\psi$ as in 6.6), 
then $P_{B,B'}=0$, $\tP_{B,B'}=0$, $\ti\L_{B,B'}=0$.}
\nl
We can find a function $\fB@>>>\ZZ$, $B\m c_B$ such that for $B,B'$ in $\fB$ 
we have $c_B=c_{B'}$ if and only if $\psi(B)=\psi(B')$. From 13.7 we deduce
$$\align&v^{c_B-c_{B'}}v^{-\k(B)-\k(B')}(B:B')\\&=
\sum_{B_1\in\fB,B_2\in\fB}
v^{c_B-c_{B_1}}P_{B_1,B}v^{c_{B_1}-c_{B_2}}\ti\L_{B_1,B_2}
v^{c_{B_2}-c_{B'}}\tP_{B_2,B'}.\endalign$$
When $B,B'$ vary, this again can be expressed as the decomposition of the
matrix $\ti\cm:=(v^{c_B-c_{B'}}v^{-\k(B)-\k(B')}(B:B'))$
(indexed by $\fB\T\fB$) as a product of three matrices $\ti\cs\ti\ct\ti\cs'$
where $\ti\cs$ (resp. $\ti\cs'$) is the matrix indexed
by $\fB\T\fB$ whose $(B,B')$-entry is 
$v^{c_B-c_{B'}}P_{B',B}$ (resp. $v^{c_B-c_{B'}}\tP_{B,B'}$) and
$\ti\ct$ is the matrix indexed by $\fB\T\fB$ whose $(B,B')$-entry is 
$v^{c_B-c_{B'}}\ti\L_{B,B'}$. 
From 7.9(a) we know that $(B:B')=0$ unless $\psi(B)=\psi(B')$. Hence
$v^{c_B-c_{B'}}(B:B')=(B:B')$ for all $B,B'$ so that $\ti\cm=\cm$. Thus we
have 
$$\ti\cs\ti\ct\ti\cs'=\cs\ct\cs'.$$
Now by 13.7(c),(e), the matrix $\ct$ (and hence the matrix $\ti\ct$) belongs
to a subgroup of $GL_N$ ($N=\sha(\fB))$ of the form $GL_{N_1}\T\do GL_{N_k}$ 
where $N_1,\do,N_k$ are the sizes of the various
subsets $\fB_\co$; moreover, by 13.7(a),(b), the matrix $\cs$ (hence the
matrix $\ti\cs$) belongs to the unipotent radical of a parabolic subgroup 
of $GL_N$ with 
Levi subgroup equal to the subgroup of $GL_N$ considered above, while the 
matrix $\cs'$ (hence the matrix $\ti\cs'$) belongs to the unipotent radical
of the opposite parabolic subgroup with the same Levi subgroup. This forces
the equalities $\ti\cs=\cs$, $\ti\ct=\ct$, $\ti\cs'=\cs'$.
Now the equality $\ti\cs=\cs$ implies $v^{c_B-c_{B'}}P_{B',B}=P_{B',B}$ for 
all $B',B$ in $\fB$. Thus, if $\psi(B)\ne\psi(B')$ (so that $c_B\ne c_{B'}$) 
we must have $P_{B',B}=0$. Similarly, from $\ti\ct=\ct$ we see that,
if $\psi(B)\ne\psi(B')$ then $\ti\L_{B,B'}=0$ and from $\ti\cs'=\cs'$ we see 
that, if $\psi(B)\ne\psi(B')$ then $\tP_{B,B'}=0$. This proves (a).

\subhead 13.9\endsubhead
We now fix $\x\in\un\fT_\et$. We define four matrices ${}^\x\cm$,
${}^\x\cs,{}^\x\ct,{}^\x\cs'$
indexed by ${}^\x\fB\T{}^\x\fB$ as follows. For $B,B'$ in ${}^\x\fB$, the $(B,B')$-entry of
${}^\x\cm$ is $v^{-\k(B)-\k(B')}(B:B')$; the $(B,B')$-entry of
${}^\x\cs$ is $P_{B',B}$; the $(B,B')$-entry of ${}^\x\ct$ is
$\ti\L_{B,B'}$; the $(B,B')$-entry of ${}^\x\cs'$ is $\tP_{B,B'}$.
Using 13.8(a) we deduce from 13.7(d) the equality of matrices 
$${}^\x\cm=({}^\x\cs)({}^\x\ct)({}^\x\cs').\tag a$$
As in the proof of 13.8(a) the last equality determines uniquely the matrices
${}^\x\cs,{}^\x\ct,{}^\x\cs'$ if the matrix ${}^\x\cm$ is known; in fact, it provides an algorithm
for computing the entries of these three matrices (and in particular for the entries
$P_{B',B}$ in terms of the entries of ${}^\x\cm$. Now under the bijection 
$\g:\BB@>\si>>{}^\x\fB$ (see 11.12(a)) the matrix ${}^\cm$ becomes a matrix indexed by
$\BB\T\BB$ whose $(b,b')$-entry is 
$v^{-\k(B)-\k(B')}(b,b')$; these entries are explicitly computable from the combinatorial
description of $(:)$ on $\VV$.
We see that

(b) {\it the quantities $P_{B',B}$ are computable by an algorithm provided that the
bijection $\g:\BB@>\si>>{}^\x\fB$ is known.}
\nl
This can be used in several examples to compute the $P_{B',B}$ explicitly.
The algorithm in (b) seems to depend on the choice of $\et$ such that $\un\et=\da$; but in
fact, by the results in 3.9, 10.22, 11.13, it does not depend on this choice.

\mpb

In the case where $m=1$, there is another algorithm to compute the quantities
$P_{B',B}$, see \cite{\CSV, Theorem 24.8}. It again displays the matrix ${}^\x\cs$ as the first
of the three factors in a matrix decomposition like (a), but with the matrix ${}^\x\cm$
being replaced by a matrix indexed by a pair of irreducible representation of a Weyl group
and with entries determined by a prescription quite different from that used in this paper.
In that case the bijection $\g:\BB@>\si>>{}^\x\fB$ is replaced by the ``generalized
Springer correspondence''. It would be interesting to understand better the connection 
between these two approaches to the quantities $P_{B',B}$.

\head 14. Odd vanishing\endhead
In this section we show that for any irreducible local system $\cl$ on a 
$G_{\un0}$-orbit in $\fg^{nil}_\da$ the cohomology sheaves of 
$\cl^\sha\in\cd(\fg_\da)$ are zero in odd degrees. (See Theorem 14.10.) 
In the case where $m\gg0$, this follows from the analogous result in the $\ZZ$-graded case in \cite{\GRA}.
In the case where $m=1$ this follows from \cite{\CSV, Theorem 24.8(a)}.
In the case where $m>1,\d=\un1$, $\cl=\bbq$ and $G$, $\fg=\op_i\fg_i$ are as in 0.3, 
this follows from \cite{\CB, Theorem 11.3} and from the known odd vanishing result
for affine Schubert varieties.

\subhead 14.1\endsubhead
We preserve the setup of 10.1, 10.2. For $\vp\in\EE'$ let $h(\vp)\in\ZZ$ be as
in 10.2. For $\boc\in\un{\ovsc\EE}$ we set $h(\boc)=h(\vp)$ where $\vp$ is any
element of $\boc\cap\EE'$; this is well defined, by 10.7(a).
For $\boc\in\un{\ovsc\EE}$ we set $T_\boc=v^{-h(\boc)}\tT_\boc\in\VV'$. Note
that $\{T_\boc;\boc\in\un{\ovsc\EE}\}$ is a $\QQ(v)$-basis of $\VV'$.
Let $f\m f^\he$ be the field involution of $\QQ(v)$ which carries $v$ to 
$-v$; this extends to a field involution of $\QQ((v))$ (denoted again by
$f\m f^\he$) given by $\sum_ic_iv^i@>>>\sum_ic_i(-v)^i$ where $c_i\in\QQ$.

Let $b\m b^\he$ be the $\QQ$-linear involution $\VV'@>>>\VV'$ such that
$(fT_\boc)^\he=f^\he T_\boc$ for any $\boc\in\un{\ovsc\EE}$ and $f\in\QQ(v)$.

\proclaim{Lemma 14.2}For any $b,b'$ in $\VV'$ we have
$(b^\he:b'{}^\he)=(b:b')^\he$. 
\endproclaim
It is enough to show that for any $\boc,\boc'$ in $\un{\ovsc\EE}$ and any
$f,f'$ in $\QQ(v)$ we have 
$$(f^\he T_\boc:f'{}^\he T_{\boc'})=(fT_\boc:f'T_{\boc'})^\he.$$
We can assume that $f=f'=1$. It is enough to show that
$$(T_\boc:T_{\boc'})\in\QQ(v^2),\tag a$$
or that
$$v^{-h(\vp_1)+h(\vp_2)}[\vp_1|\vp_2]\in\QQ(v^2)$$
for any $\vp_1,\vp_2$ in $\EE'$ (notation of 10.12(a)) or that
$$v^{-h(\vp_1)+h(\vp_2)}\sum_{w\in\cw}v^{\t(\vp_2,w\vp_1)}\in\QQ(v^2).\tag b$$
From 10.8(b) we see that for $\vp,\vp'$ in $\EE'$ we have 
$$\t(\vp,\vp')=h(\vp)+h(\vp')\mod2.\tag b$$
From 10.11(b) we see that for $\vp\in\EE'$, $w\in\cw$ and $N\in\ZZ$ we have
$$\dim{}^\e\fp_N^{\un{w\vp}}=\dim{}^\e\fp_N^{\un{\vp}},\text{ hence }
\dim{}^\e\fu_N^{\un{w\vp}}=\dim{}^\e\fu_N^{\un{\vp}},$$
so that
$$h(w\vp)=h(\vp).$$
Using this and (b) we see that for $\vp,\vp'$ in $\EE'$ we have
$$\sum_{w\in\cw}v^{\t(\vp',w\vp)}\in v^{h(\vp')+h(\vp)}\NN[v^2,v^{-2}].$$
Thus, (b) holds. The lemma is proved.

\subhead 14.3\endsubhead
From 14.2 we deduce that ${}^\he:\VV'@>>>\VV'$ carries $\fR_l=\fR_r$ onto 
$\fR_l=\fR_r$; hence
it induces a $\QQ$-linear involution $\VV@>>>\VV$ (denoted again by ${}^\he$).
From 14.2 we deduce:

(a) {\it For any $b,b'$ in $\VV$ we have $(b^\he:b'{}^\he)=(b:b')^\he$.}

\subhead 14.4\endsubhead
Let $\boc\in\un{\ovsc\EE},\vp\in\EE''\cap\boc$ and let 
$a_\vp=\sum_jv^{-2s_j}v^{d(\vp)}\in\ca$ be as in 11.1. Note that 
$$(a_\vp)^\he=\sum_jv^{-2s_j}(-v)^{d(\vp)}=(-1)^{d(\vp)}a_\vp.$$
Hence we have
$$(a_\vp\i\tT_\boc)^\he=(a_\vp\i v^{h(\boc)}T_\boc)^\he=
(-1)^{d(\vp)+h(\boc)}a_\vp\i v^{h(\boc)}T_\boc=
=(-1)^{d(\vp)+h(\boc)}a_\vp\i\tT_\boc.$$
It follows that 

(a) {\it ${}^\he:\VV'@>>>\VV'$ carries $\VV'_\ca$ onto $\VV'_\ca$ and
${}^\he:\VV@>>>\VV$ carries $\VV_\ca$ onto $\VV_\ca$.}

\subhead 14.5\endsubhead
We show:

a) {\it For any $b\in\VV'$ we have $\ov{b^\he}=(\bar b)^\he$. Hence for any 
$b\in\VV$ we have $\ov{b^\he}=(\bar b)^\he$.}
\nl
We can assume that $f=fT_\boc$ where $\boc\in\un{\ovsc\EE}$ and $f\in\QQ(v)$.
Note that $\ov{f^\he}=(\baf)^\he$. Hence we can assume that $f=1$. We have
$$\ov{T_\boc}=\ov{v^{-h(\boc)}\tT_\boc}=v^{h(\boc)}\tT_\boc
=v^{2h(\boc)}T_\boc,$$  
hence
$$(\ov{T_\boc})^\he=v^{2h(\boc)}T_\boc.$$
We have $\ov{T_\boc^\he}=\ov{T_\boc}=v^{2h(\boc)}T_\boc$.
This proves (a).

\subhead 14.6\endsubhead
We show:

(a) {\it $b\m b^\he$ is a bijection $\BB'@>\si>>\BB'$.}
\nl
Let $b\in\BB'$. From $b\in\VV_\ca$ we see using 14.4(a) that $b^\he\in\VV_\ca$.
From $\bar b=b$ we see using 14.5(a) that $\ov{b^\he}=(\bar b)^\he=b^\he$. 
From $(b:b)\in\QQ(v)\cap(1+v\ZZ[[v]])$ we see using 14.3(a) that
$$(b^\he:b'{}^\he)=(b:b')^\he\in(\QQ(v)\cap(1+v\ZZ[[v]]))^\he=
\QQ(v)\cap(1+v\ZZ[[v]])\sub1+v\ZZ[[v]].$$
Using this and the definitions, we see that $b^\he\in\BB'$. Thus the map 
$b\m b^\he$, $\BB'@>>>\BB'$ is well defined. Since this map has square $1$, it
is a bijection. This proves (a).

\mpb
From (a) we deduce:

(b) {\it If $b\in\BB$ then $b^\he=s_b\ti b$ for a well defined 
$s_b\in\{1,-1\}$ and a well defined $\ti b\in\BB$.}
\nl
The following result makes (b) more precise.

\proclaim{Lemma 14.7}Let $\co$ be a $G_{\un0}$-orbit in $\fg_\da^{nil}$ and let
$b\in\BB_\co$. We have $b^\he=(-1)^{\dim\co}b$.
\endproclaim
We argue by induction on $\dim\co$; we can assume that the result holds
when $\co$ is replaced by an orbit of dimension $<\dim\co$ (if any).
We have $\g(b)=\cl^\sha[\dim\co]$ where 
$(\co,\cl)\in\ci(\fg_\da)$. Let $x\in\co$; we associate to $x$ an $\e$-spiral 
$\fp^\ph_*$ and a splitting $\tfl^\ph_*$ of it as in 2.9. Recall that
$\ovsc\tfl^\ph_\et\sub\co$ and that $\cl_1:=\cl|_{\ovsc\tfl^\ph_\et}$ is an 
irreducible $\tL^\ph_0$-equivariant local system on $\ovsc\tfl^\ph_\et$; thus,
$\cl_1^\sha\in\cd(\tfl^\ph_\et)$ is defined. Let 
$I={}^{\e}\Ind_{\fp_\et^{\phi}}^{\fg_\da}(\cl_1^\sha)\in\cq(\fg_\da)$;
we have clearly $I\in{}^\x\cq(\fg_\da)$. 
By 2.9(b),(c), in ${}^\x\ck(\fg_\da)$ we have 
$$I=\cl^\sha+\sum_{\co';\dim\co'<\dim\co}\sum_{b'\in\BB_{\co'}}f_{b'}\g(b')$$
where $f_{b'}\in\ca$. Define $I'\in\VV_\ca$ by $\g(I')=I$. We have
$$I'=v^{-\dim\co}b+\sum_{\co';\dim\co'<\dim\co}\sum_{b'\in\BB_{\co'}}f_{b'}b'.
\tag a$$
By \cite{\GRA, 17.2, 17.3}, in $\ck(\tfl^\ph_\et)$ we have
$$h\cl_1^\sha=\sum_{j\in J}h_j\ind_{\fq(j)_\et}^{\tfl^\ph_\et}(C(j)),\tag b$$
where $h=h^\he\in\ca-\{0\}$, $J$ is a finite set and for $j\in J$, $\fq(j)$ is
a parabolic subalgebra of $\tfl^\ph_\et$ with Levi subalgebra $\fm(j)$ such 
that $\fq(j),\fm(j)$ are compatible with the
$\ZZ$-grading of $\tfl^\ph$, $C(j)\in\cq(\fm(j)_\et)$ is a cuspidal
perverse sheaf and $h_j\in\ca$. 
Moreover, since $\cl_{1}^{\sha}$ belongs to the block of $\cq(\tfl_{\et}^{\phi})$ given by $\x$, we can assume that $\fm(j)=\fm$ and 
$C(j)=\tC$ for all $j$. Thus (b) can be written in the form
$$h\cl_1^\sha=F_0+F_1,\tag c$$
where
$$F_0=\sum_{j\in J}h'_j\ind_{\fq(j)_\et}^{\tfl^\ph_\et}(\tC[-\dim\fm_\et])),$$
$$F_1=\sum_{j\in J}h''_j\ind_{\fq(j)_\et}^{\tfl^\ph_\et}(\tC[-\dim\fm_\et])),$$
$$h'_j+h''_j=h_jv^{\dim\fm_\et},h'_j\in\ZZ[v^2,v^{-2}],
h''_j\in v\ZZ[v^2,v^{-2}].$$ 
Let $\ck(\tfl_\et^{\phi})^{ev}$ be the $\ZZ[v^2,v^{-2}]$-submodule of 
$\ck(\tfl_\et^{\phi})^{ev}$ with basis $\cl'{}^\sha$ for various $(\co',\cl')$ in 
$\ci(\tfl_\et^{\phi})$. By \cite{\GRA, 21.1(c)}, for any $j\in J$, we have 
$$\ind_{\fq(j)_\et}^{\tfl^\ph_\et}(\tC[-\dim\fm_\et])\in\ck(\tfl_{\et}^{\phi})^{ev}.$$
Hence $F_0\in\ck(\tfl_{\et}^{\phi})^{ev}, F_1\in v\ck(\tfl_{\et}^{\phi})^{ev}$. We have clearly
$h\cl_1^\sha\in\ck(\tfl_{\et}^{\phi})^{ev}$. Since
$v\ck(\tfl_{\et}^{\phi})^{ev}\cap\ck(\tfl_{\et}^{\phi})^{ev}=0$, we deduce from (c) that
$h\cl_1^\sha=F_0$, that is
$$h\cl_1^\sha=\sum_{j\in J}h'_j\ind_{\fq(j)_\et}^{\tfl^\ph_\et}
(\tC[-\dim\fm_\et]))$$
where $h'_j\in\ZZ[v^2,v^{-2}]$. Applying ${}^\e\Ind_{\fp_\et}^{\fg_\da}$ and 
using the transitivity property 4.2 (and also 10.4(a)) we obtain
$$hI=\sum_{j\in J}h'_jI_{\vp_j}$$
where $\vp_j\in\EE'$ for $j\in J$. Since $\g\i(\sum_{j\in J}h'_jI_{\vp_j})$ is
fixed by ${}^\he$, it follows that $hI'$ is fixed by ${}^\he$. Since $h\ne0$
and $h^\he=h$, we deduce that $I'{}^\he=I'$. Using (a) and the induction 
hypothesis, we see that
$$\align&(v^{-\dim\co}b)^\he+
\sum_{\co';\dim\co'<\dim\co}\sum_{b'\in\BB_{\co'}}f_{b'}^\he(-1)^{\dim\co'}b'
\\&=v^{-\dim\co}b+\sum_{\co';\dim\co'<\dim\co}\sum_{b'\in\BB_{\co'}}f_{b'}b'.
\endalign$$
Using this and 14.6(b), we deduce that

$(-v)^{-\dim\co}s_b\ti b-v^{-\dim\co}b$ is a $\ca$-linear combination of 
elements in 

$\cup_{\co';\dim\co'<\dim\co}\BB_{\co'}$.
\nl
Since $b\in\BB_\co$ and $\ti b\in\BB$ this forces $\ti b=b$ and 
$s_b=(-1)^{\dim\co}$. This completes the proof of the lemma.

\subhead 14.8\endsubhead
We show:

(a) {\it Let $\co,\co'$ be $G_{\un0}$-orbits in $\fg_\da^{nil}$ and let $B,B'$ 
in $\fB$ be such that the support of $B$ (resp. $B'$) is the closure of $\co$ 
(resp. $\co'$). We have $(B:B')^\he=(-1)^{\dim\co+\dim\co'}(B:B')$.}
\nl
If $\psi(B)\ne\psi(B')$ then $(B:B')=0$ by 7.9(a) and (a) holds. Assume now 
that
$\psi(B)=\psi(B')$. We can assume that $\psi(B)=\psi(B')=\x$. It is enough to
prove that, if $b\in\BB_\co$, $b'\in\BB_{\co'}$, then
$(b:b')^\he=(-1)^{\dim\co+\dim\co'}(b:b')$. Using 14.2 and 14.7, we have
$$(b:b')^\he=(b^\he:b'{}^\he)=((-1)^{\dim\co}b,(-1)^{\dim\co'}b')$$
and (a) is proved.

\subhead 14.9\endsubhead
We show: 

(a) {\it For any $B,B'$ in $\fB$ we have $P_{B,B'}\in\NN[v^{-2}]$,
$\tP_{B,B'}\in\NN[v^{-2}]$, $\ti\L_{B,B'}\in\QQ((v^2))$ (notation of 13.7).}
\nl
The proof is similar to that of 13.8(a). With the notation in 14.8(a) we apply
${}^\he$ to
$$v^{-\dim\co-\dim\co'}(B:B')=
\sum_{B_1\in\fB,B_2\in\fB}P_{B_1,B}\ti\L_{B_1,B_2}\tP_{B_2,B'}$$ 
(see 13.7(d)); using 14.8 we obtain
$$v^{-\dim\co-\dim\co'}(B:B')=
\sum_{B_1\in\fB,B_2\in\fB}P_{B_1,B}^\he\ti\L_{B_1,B_2}^\he\tP_{B_2,B'}^\he.$$
When $B,B'$ vary, this can be expressed as the decomposition of the
matrix $\cm=(v^{c_B-c_{B'}}v^{-\k(B)-\k(B')}(B:B'))$
(indexed by $\fB\T\fB$) as a product of three matrices $\cs^{\he}\ct^{\he}{\cs'}^{\he}$ 
where $\cs^{\he}$ (resp. ${\cs'}^{\he}$) is the matrix indexed
by $\fB\T\fB$ whose $(B,B')$-entry is $P_{B',B}^\he$ (resp. $\tP_{B,B'}^\he$) 
and $\ct^{\he}$ is the matrix indexed by $\fB\T\fB$ whose $(B,B')$-entry is 
$\ti\L_{B,B'}^\he$. Recall from 13.7 that we have also $\cm=\cs\ct\cs'$ 
(notation of 13.7). Thus we have 
$$\cs^{\he}\ct^{\he}{\cs'}^{\he}=\cs\ct\cs'.$$
Now by 13.7(c),(e), the matrix $\ct$ (and hence the matrix $\ct^{\he}$) belongs
to a subgroup of $GL_N$ ($N=\sha(\fB))$ of the form $GL_{N_1}\T\do GL_{N_k}$ 
where $N_1,\do,N_k$ are the sizes of the various
subsets $\fB_\co$; moreover, by 13.7(a),(b), the matrix $\cs$ (hence the
matrix $\cs^{\he}$) belongs to the unipotent radical of a parabolic subgroup of
$GL_N$ with Levi subgroup equal to the subgroup of $GL_N$ considered above, 
while the matrix $\cs'$ (hence the matrix ${\cs'}^{\he}$) belongs to the unipotent
radical of the opposite parabolic subgroup with the same Levi subgroup. This 
forces the equalities $\cs^{\he}=\cs$, $\ct^{\he}=\ct$, ${\cs'}^{\he}=\cs'$. Now the 
equality $\cs^{\he}=\cs$ implies $P_{B',B}^\he=P_{B',B}$ for all $B',B$ in 
$\fB$. Similarly, from $\ct^{\he}=\ct$ we see that 
$\ti\L_{B,B'}^\he=\ti\L_{B,B'}$ for all $B,B'$ in $\fB$ and from 
${\cs'}^{\he}=\cs'$ we see that $\tP_{B,B'}^\he=\tP_{B,B'}$. This proves (a).

\proclaim{Theorem 14.10} (a) Let $(\co,\cl)\in\ci(\fg_\da)$ and let 
$A=\cl^\sha\in\cq(\fg_\da)$.
We have $\ch^aA=0$ for any odd integer $a$.

(b) Let $\vp\in\EE'$. We have $\ch^a(I_\vp)=0$ for any odd integer $a$.

(c) Let $(\fp_*,L,P_0,\fl,\fl_*,\fu_*)$ be as in 4.1(a) with $\e=\doet$, and let
$(\co',\cl')\in\ci(\fl_\et)$ (see 1.2). We form $A'=\cl'{}^\sha\in\cq(\fl_\et)$. 
We have $\ch^a({}^\e\Ind_{\fp_\et}^{\fg_\da}(A'))=0$ for any odd integer $a$.

(d) In the setup of (c), ${}^\e\Ind_{\fp_\et}^{\fg_\da}(A')$ is a direct sum of
complexes of the form $\cl^\sha[s]$ for various $(\co,\cl)\in\ci(\fg_\da)$ and
various even integers $s$.
\endproclaim
(a) follows from 14.9(a).

We prove (b). Let $\boc\in\un{\ovsc\EE}$ be such that $\vp\in\boc$.
In $\VV$ we have $T_\boc=\op_{b\in\BB}f_bv^{-\k(b)}b$ where $f_b\in\ca$. 
Applying ${}^\he$ and using $T_\boc^\he=T_\boc$ 
and $(v^{-\k(b)}b)^\he=v^{-\k(b)}b$ for $b\in\BB$ (see 14.7) we see that
$$\op_{b\in\BB}f_b^\he v^{-\k(b)}b=\op_{b\in\BB}f_bv^{-\k(b)}b.$$
Hence $f_b^\he=f_b$ that is, $f_b\in\ZZ[v^2,v^{-2}]$. Thus,

(e) {\it $I_\vp$ is isomorphic to a direct sum of complexes of the form
$B[-\k(B)][2s]$ with $B\in\fB$ and $s\in\ZZ$.}
\nl
Hence $\ch^aI_\vp$ is isomorphic to
a direct sum of sheaves of the form $\ch^{a+2s}(B[-\k(B)])$ with $B\in\fB$. 
Hence the desired result follows from (a).

We prove (c). By 1.5(a) we can find $\fq_*$, 
$(\hfp_*,\hL,\hP_0,\hfl,\hfl_*,\hfu_*)$ as in 4.2 with $\e=\doet$ and a cuspidal
perverse sheaf $C$ in $\cq(\hfl_\et)$ such that $A'[s']$ is a direct summand of
$\ind_{\fq_\et}^{\fl_\et}(C[-\dim\fl_\et])$ for some $s'\in\ZZ$ hence, by 4.2(a),
${}^\e\Ind_{\fp_\et}^{\fg_\da}(A')[s']$ is a direct summand of
${}^\e\Ind_{\hfp_\et}^{\fg_\da}(C[-\dim\fl_\et])$; moreover, by 
\cite{\GRA, 21.1(c)}, we have $s'=2s''$ for some $s''\in\ZZ$. Hence
$\ch^a({}^\e\Ind_{\fp_\et}^{\fg_\da}(A'))$ is a direct summand of
$\ch^{a-2s''}({}^\e\Ind_{\hfp_\et}^{\fg_\da}(C[-\dim\fl_\et]))$. We can assume
that $\hfp_\et$ is an $\e$-spiral with splitting $\fm_*$ (in 10.1) and that
$C=\tC$ (in 10.1). Then
${}^\e\Ind_{\hfp_\et}^{\fg_\da}(C[-\dim\fl_\et])=I_\vp$ for some $\vp\in\EE'$ 
and $\ch^a({}^\e\Ind_{\fp_\et}^{\fg_\da}(A'))$ is a direct summand of
$\ch^{a-2s''}I_\vp$. The desired result follows from (b).

We prove (d). As in the proof of (c), ${}^\e\Ind_{\fp_\et}^{\fg_\da}(A')$ is a 
direct summand of $I_\vp[-2s'']$ for some $\vp\in\EE'$ and $s''\in\ZZ$. This,
together with (e) gives the desired result. The theorem is proved.


\widestnumber\key{ABC}
\Refs
\ref\key{\KLL}\by D.Kazhdan and G.Lusztig\paper Schubert varieties and 
Poincar\'e duality\jour Proc. Symp. Pure Math.\vol36\paperinfo Amer. Math. 
Soc.\yr1980\pages185-203\endref
\ref\key{\GRA}\by G.Lusztig\paper Study of perverse sheaves arising from graded
Lie algebras\jour Adv.Math.\vol112\yr1995\pages147-217\endref
\ref\key{\CUS}\by G.Lusztig\paper Cuspidal local systems and graded Hecke
algebras\jour Publ.Math. IHES\vol67\yr1988\pages145-202\endref
\ref\key{\CSV}\by G.Lusztig\paper Character sheaves, V\jour Adv. in Math.\vol64\yr1986
\pages103-155\endref
\ref\key{\CB}\by G.Lusztig\paper Canonical bases arising from quantized enveloping
algebras \jour Jour. Amer. Math. Soc.\vol3\yr1990\pages447-498\endref
\ref\key{\LY}\by G.Lusztig and Z.Yun\paper $\ZZ/m$-graded Lie algebras and 
perverse sheaves, I\jour arxiv:1602.05244\endref
\endRefs
\enddocument